\documentclass[a4paper,11pt, twoside]{article}

\usepackage{amsfonts,amsthm,amsmath,amsfonts,latexsym,amssymb}

\usepackage{verbatim}
\usepackage{mathtools}
\usepackage{enumerate}
\usepackage{mathrsfs}
\usepackage{dsfont}
\usepackage{bbm}
\usepackage{vmargin}
\usepackage[]{appendix}
\usepackage{enumitem}
\usepackage{soul}
\usepackage[normalem]{ulem}
\usepackage{cancel}
\usepackage{xcolor}
\usepackage{authblk} 
\definecolor{naranja}{rgb}{1,0.5,0}

\newcommand{\tachar}[1]{%
  \ifmmode
    {\color{naranja}\cancel{#1}}%
  \else
    {\color{naranja}\sout{#1}}%
  \fi
}
\usepackage[utf8]{inputenc}
\usepackage[T1]{fontenc}

\usepackage{graphicx,upref,calligra,float}
\usepackage[rflt]{floatflt}
\usepackage[svgnames]{xcolor}

\usepackage[colorlinks=true,linkcolor=azul,citecolor=azul]{hyperref}
\usepackage{fancyhdr}
\usepackage{lastpage}
\usepackage{etoolbox}
\usepackage{tocloft}  \renewcommand\numberline[1]{}

\usepackage[backend=bibtex,doi=false,url=false,eprint=false,isbn=false,maxbibnames=7,sorting=nty,giveninits=true]{biblatex}
\newtheorem{teor}{Theorem}[section]

\newtheorem{fed}[teor]{Definition}
\theoremstyle{definition}
\newtheorem{teo}[teor]{Theorem}
\newtheorem*{teo*}{Theorem}
\newtheorem*{obs*}{Remark}
\newtheorem{lem}[teor]{Lemma}
\newtheorem{pro}[teor]{Proposition}
\newtheorem{cor}[teor]{Corollary}

\newtheorem{obs}[teor]{Remark}
\newtheorem{prop}[teor]{Proposition}

\theoremstyle{definition}

\newtheorem{rem}[teor]{Remark}

\newtheorem{exa}[teor]{Example}

\definecolor{azul}{rgb}{0.1,0.6,0.86}
\definecolor{naranja}{RGB}{249,153,96}
\definecolor{titleblue}{rgb}{0.13,0.49,0.69}
\definecolor{mylred}{rgb}{0.85,0.24,0.2}
\definecolor{myblue}{rgb}{0,0.33,0.55}

\def\bdem{\begin{proof}}
\def\edem{\end{proof}}
\def\ds{\displaystyle}

\def\eps{\varepsilon}
\def\fii{\varphi }
\def\la{\lambda}

\def\w{\omega}

\def\a{\alpha}

\def\N{\mathbb{N}}

\def\R{\mathbb{R}}

\def\C{\mathbb{C}}

\def\E{\mathbb{E}}

\def\cC{\mathcal{C}}
\def\cD{\mathcal{D}}

\def\cP{\mathcal{P}}

\newcommand{\ce}{\mathbb{C}}
\newcommand{\e}{\mathbb{E}}

\newcommand{\efe}{\mathcal{F}}
\newcommand{\pe}{\mathbb{P}}

\DeclareMathOperator{\var}{Var}
\DeclareMathOperator{\cov}{Cov}
\DeclareMathOperator{\re}{Re}

\newcommand{\deriv}{\mathrm{d}}

\newcommand{\schw}{\mathcal{S}(\R^d)}
\newcommand{\schwd}{\mathcal{S}'(\R^d)}

\title{Stability of crystalline measures under random perturbations}
\author[1,2]{Jorge Antezana\thanks{Corresponding author. E-mail \texttt{jorge.antezana@ub.edu}. Supported by Grant PID2024-160033NB-I00 funded by MICIU/AEIMCIN/AEI/10.13039/501100011033 and by /10.13039/501100011033FEDER, UE.}}
\author[3,4]{Luciano Scazzola \thanks{E-mail \texttt{lscazzola@mate.unlp.edu.ar} . Supported by National University of La Plata Grant 11X974, CONICET-PIP 11220210100954CO }}

\affil[1]{Departament de Matemàtiques i Informàtica, Universidad de Barcelona, Barcelona, Spain.}
\affil[2]{Centre de Recerca Matemàtica, Barcelona, Spain.}
\affil[3]{Universidad Nacional de La Plata, La Plata, Argentina.}
\affil[4]{Instituto Argentina de Matemática ``Alberto P. Calder\'on'', CONICET, Argentina.}

\date{}

\begin{document}

\maketitle

\begin{abstract}


In this paper, we study random perturbations and the diffraction properties of translation-bounded crystalline measures. Under a non-vanishing assumption on the known perturbation factor, we prove that, almost surely, any realization of the perturbed measure preserves sufficient structural information to allow the full recovery of the original measure. Previous works have only considered crystalline measures with uniformly discrete supports and equal-mass atoms. Our general setting includes many important crystalline measures that were excluded until now. To treat this broader setting, we adopt a different framework. In this work, we model the perturbation as a random measure rather than as a random set. In particular, this leads us to consider not only perturbations of the support, but also perturbations of the atoms' masses. On the other hand, using elementary techniques of Fourier analysis, we prove that any translation-bounded crystalline measure has a diffraction measure, and we explicitly compute its diffraction.

\end{abstract}

\tableofcontents

\section{Introduction}

The discovery of quasicrystals by Shechtman in 1982 revealed the existence of solids exhibiting sharp Bragg diffraction patterns with rotational symmetries forbidden by classical crystallography. This led to the development of mathematical models for aperiodic structures that describe these new materials. Among them, \emph{model sets} (or \emph{Meyer sets}), introduced by Meyer almost a decade earlier, play a central role in the mathematical description of quasicrystals.

\medskip

A model set can be regarded as a measure by placing equal-mass atoms at each of its points, thereby defining a purely atomic measure with uniformly separated support. In particular, it induces a tempered distribution. However, the distributional Fourier transform of such a measure is typically supported on a (countable) dense set. On the other hand, Lev and Olevskii \cite{MR3338010} studied purely atomic measures whose support and spectral support are both \emph{uniformly discrete}, and showed that such configurations must essentially arise from generalized Dirac combs. This established a rigidity result: the coexistence of uniform separation in both physical and frequency space forces the measure to be closely related to a lattice structure.

\medskip

This naturally raises the question of whether both the support of a measure and that of its Fourier transform can merely be \emph{locally finite}, without being uniformly discrete. This led to the following definition.

\medskip

\begin{fed}[Meyer \cite{MR3482845}]
A purely atomic measure $ \mu $ defined on the Borel sets of $ \mathbb{R}^d$ is called a \textbf{crystalline measure} if:
\begin{enumerate}
    \item[i)] The support $ X $ of $ \mu $ is a locally finite set.
    \item[ii)] The measure $ \mu $ induces a tempered distribution.
    \item[iii)] The distributional Fourier transform $ \widehat{\mu} $ of $ \mu $ is also a purely atomic measure supported on a locally finite set $ X^* $.
\end{enumerate}
\end{fed}

\medskip

Note that a crystalline measure induces a \emph{generalized Poisson formula} of the form:
\begin{equation}\label{poisson}
    \sum_{p \in X} \mu(p) \, \widehat{f}(p) = \sum_{\lambda \in X^*} A(\lambda) \, f(\lambda),
\end{equation}
which holds for any function $f$ in the Schwartz class, which we will denote by $\schw $. Conversely, any such identity, with both $ X $ and $ X^* $ locally finite, gives rise to a crystalline measure.

\medskip

Besides their connection with quasicrystals, this concept has connections with a wide range of deep problems in analysis. In number theory, for instance, Kahane and Mandelbrot \cite{MR106208} proved a relationship between crystalline measures and Dirichlet series satisfying the same functional equation as the Riemann Zeta function (see also \cite{MR4939386}). These measures are also related to the spectral theory of quantum graphs, the distribution of zeros of Dirichlet polynomials, uniqueness pairs in harmonic analysis and formulations of the uncertainty principle, and interpolation formulae in the Schwartz class (see for instance \cite{MR4697523}, \cite{MR4206541}, \cite{MR4294062}, \cite{MR3949027}, respectively).

\medskip

All the aforementioned connections allow us to understand why Dyson \cite{MR2483565} highlights the importance of understanding the structural properties of such measures in his well-known essay \emph{Birds and Frogs}. However, given these rich and varied connections, it is clear that the internal structure of crystalline measures is highly intricate. In recent years, several lines of research have emerged intending to uncover and analyze this complexity (\cite{MR4664429}, \cite{MR4841781}, \cite{Fa24}, \cite{ArXiv-2312.11185}, \cite{MR4678143}, among others). The present work fits within this broader program. Ramos and Sousa \cite{MR4678143} studied the stability of function reconstruction from perturbed sampling nodes using functional-analytic techniques. Here, we investigate how crystalline measures, taking both atomic locations and masses into account, behave under random perturbations and study their diffraction properties using diffraction-theoretic methods. Our work is more directly motivated by recent results on the perturbation of structured distributions of points \cite{MR4526342,MR4406130}.

\medskip

Yakir \cite{MR4406130} considered a Dirac comb associated with a lattice in which each atom is randomly displaced using independent and identically distributed vectors. Relying on the \emph{group structure} of the lattice and \emph{the classical Poisson summation formula}, he showed that almost surely it is possible to recover the original configuration from almost every realization. Of course, some natural moment conditions on the perturbation are imposed.

\medskip

This result leads to the question of whether similar recovery phenomena hold beyond the periodic setting. A first step in this direction was taken by Petrache and Viera~\cite{MR4526342}, who extended Yakir’s framework to atomic measures supported on uniformly discrete sets. In the absence of a group structure, they employed techniques from \emph{diffraction theory} to control the averaged behavior of the perturbed configuration and succeeded in establishing a recovery result in this more general context. However, the assumption that the \emph{support} of the measure is uniformly discrete remains a strong geometric constraint, which excludes many important examples of crystalline measures such as those associated with the zeros of Dirichlet polynomials \cite{MR4664429, MR4841781, MR4206541}, spectral data of differential operators on quantum graphs \cite{MR4129870}, the families of measures independently constructed by Kolountzakis \cite{MR3448912} and Meyer \cite{MR3482845}, and positive crystalline measures with masses bounded from below, whose supports are typically finite unions of uniformly separated sets.

\subsection{Framework}\label{framework}

Recall that a \emph{Radon measure} $\mu$ on $\R^d$ (henceforth simply called a \textit{measure}) is a linear functional on the vector space of compactly supported continuous functions $C_c(\mathbb{R}^d)$ such that for every compact set $K \subseteq \mathbb{R}^d$ there is a positive constant 
$M_K$ satisfying that
\[
|\mu(f)| \leq M_K\|f\|_\infty \,,
\]
for every $f \in C_c(\mathbb{R}^d)$ supported on $K$, and where $\|\cdot\|_\infty$ denotes the 
supremum norm on $C_c(\mathbb{R}^d)$. We say $\mu$ is a \textit{positive measure} if $\mu(f)\geq 0$ for every non-negative $f\in C_c(\R^d)$. The \textit{total variation} $|\mu|$ of a measure $\mu$ is the positive measure characterized, for every non-negative $f\in C_c(\R^d)$, by
\[
 |\mu|(f) = \sup\{|\mu(g)| : g\in C_c(\R^d), |g|\leq f\} 
\]
(see \cite{MR3791850}). By
the Riesz Representation Theorem, there is an equivalence between this definition of measure and
the classical measure-theoretic concept of regular Radon measure. 

\medskip

Throughout this paper, $\|\cdot\|$ will denote the euclidean norm in $\R^d$. Given $x_0\in\R^{d}$ and $R>0$, let
$$
B_R(x_0):=\{x\in\R^d: \lVert x-x_0\rVert\leq R\},
$$
and let $B_R := B_R(0)$.  A (Radon) measure $\mu$ is called \textit{translation-bounded} if
\[
\sup_{x \in \mathbb{R}^d} |\mu|(B_1(x)) < \infty \,.
\]

In the present work, we will study random perturbations and diffraction properties of crystalline measures. The natural framework for investigating both problems is the one provided by the translation-bounded measures. This is a much weaker assumption than those considered in previous works, where the support was assumed uniformly separated, and all atoms had equal mass. The hypothesis of translation boundedness not only allows us to study more general supports, but also measures with atoms of different masses. These masses can be real or complex.  This broader setting also includes all the previously mentioned examples, which had fallen outside the scope of existing results. There are some crystalline measures that are not translation-bounded, many of them related to sampling formulas in the Schwartz class \cite{MR3949027, MR4866154}. However, the assumption of translation boundedness is more natural in the context of diffraction theory~\cite{MR3136260,MR1328260}, and any positive crystalline measure has to be translation-bounded \cite[Prop. 3.3]{MR1328260}. Hence, we have decided to keep our attention on this class of measures. 

\textbf{Notation:} To simplify the statements, for any set $Y\subseteq\R^d$, we write $Y_R:=Y\cap B_R$.

\subsection{Main results}

\subsubsection*{Random perturbations}

Our first main theorem reads as follows.

\begin{teor}\label{principal1}
 Let $\mu$ be a \emph{translation-bounded} crystalline measure  on $\R^d$ with
 \begin{equation}\label{eq:N}\tag*{\textbf{(N)}}
\mu=\sum_{p\in X}\mu(p)\delta_p,
\qquad
\widehat{\mu}=\sum_{\lambda\in X^*}A(\lambda)\delta_\lambda \,.
\end{equation}
 Let $ \{\xi_{p}\}_{p\in X}$ be an i.i.d. family of random vectors in $\R^d$, all distributed as a random vector $\xi$ satisfying 
\[
  \E(\|\xi\|^{\max\{2,d\}})<\infty. 
\]
Set
\[
    w(\lambda):=\E\bigl(e^{-2\pi i\langle \xi,\lambda\rangle}\bigr),
    \qquad \lambda\in\R^d,
\]
and define the random measure
$$
\widetilde{\mu} :=\sum_{p\in X} \mu(p) \delta_{p + \xi_p}.
$$
Then, almost surely, for every $\la\in\R^d$ we have that
\begin{equation}\label{mainresult}
       \lim_{R\to\infty} \frac{1}{m_d(B_R)} \int_{B_R}e^{-2\pi i \langle x,\la \rangle} \,d\widetilde\mu(x) =  \left\{ \begin{array}{lcc}
             w(\la) A(\lambda)  &   if  & \lambda\in X^* \\
             \\ 0 &  if  & \lambda\notin X^*
             \end{array}.
   \right.
    \end{equation}
\end{teor}

\begin{obs*}
The identical distribution assumption in Theorem~\ref{principal1} is not essential.
In fact, the proof only requires the following more general assumptions on the random vectors:
\begin{enumerate}[
    label=\textbf{(A\arabic*)},
    ref=A\arabic*
]
    \item\label{eq:A1} The random vectors $\{\xi_p\}_{p\in X}$ are independent and there exists a function $w:\R^d\to\C$ such that, for all $\la\in\R^d$,
\[\limsup_{\|p\|\to\infty}
\left|
\E\left(e^{-2\pi i\langle\xi_p,\lambda\rangle}\right)-w(\lambda)
\right|=0 \,. \] 
    \item\label{eq:A2} There exists a random vector $\xi$, with
$\E(\|\xi\|^{\max\{2,d\}})<\infty$, such that it controls the decay of the tails of the random vectors $\xi_p$ uniformly, more precisely, for every $p\in X$ and every $\eta\geq 0$
\[
\pe(\|\xi_p\|>\eta)\le\pe(\|\xi\|>\eta) \,.
\]
\end{enumerate}
Thus, Theorem~\ref{principal1} remains valid under assumptions {\hypersetup{linkcolor=azul}%
\textnormal{\textbf{(\ref{eq:A1})}}}-{\hypersetup{linkcolor=azul}%
\textnormal{\textbf{(\ref{eq:A2})}}}, and the proof below is carried out in this more general setting.
\end{obs*}

\bigskip

As previously established by Yakir for lattices and by Petrache and Viera in their framework, the perturbation introduces an unstructured component that can be decoupled from the original measure via the Fourier transform and removed with appropriate averaging. One of the main consequences of Theorem \ref{principal1} is that, after perturbing a crystalline measure, an echo of the original crystalline measure remains in almost every realization of the perturbed measure. 
More precisely

\begin{cor}

    Under the assumptions of Theorem \ref{principal1}, let
$$
U:=\{\lambda\in\mathbb R^d:w(\lambda)\neq0\}.
$$
If \(w\) is known, then a single realization of \(\widetilde{\mu}\) almost surely determines \(\widehat{\mu}|_U\). In particular, if \(w\) is nowhere zero, it determines \(\widehat{\mu}\) and hence \(\mu\). Moreover, under this non-vanishing assumption, the support \(X^*\) can be recovered without knowing the values of \(w\). At frequencies where \(w\) vanishes, the limits in Theorem \ref{principal1} are zero regardless of \(A(\lambda)\), so they do not determine the corresponding coefficients without additional information.
\end{cor}

\medskip

At first glance, since in this theorem we are only perturbing the support of the measure, the reader may think that we are still paying attention to a random set instead of the whole measure. However, as soon as the reader goes into the proofs, he realizes that to overcome many difficulties produced by the fact that the support is not uniformly separated, the whole measure has to be taken into account. Therefore, although the general scheme of proof is the same as that followed in \cite{MR4526342,MR4406130}, new ideas and techniques are required since the problem has moved from the setting of random sets to the setting of random measures. In some cases, the new techniques allow us to slightly improve the previous results. For instance, in those works, the moment condition imposed on the random vectors $\xi_p$ is  $\E(\|\xi_p\|^{d+\eps})<\infty$. We show that the \(d\)-th moment condition suffices for the boundary-crossing estimate in Lemma \ref{adentroafuera}. However, our random-part argument requires a finite second moment, so Theorem~\ref{principal1} assumes \(\mathbb E(\|\xi\|^{\max\{2,d\}})<\infty\). In particular, this yields a slightly weaker condition when $d\geq 2$.

\bigskip

Another reason why in Theorem~\ref{principal1} we only perturb the locations of the points in the support is that the notation is simpler and the arguments are easier to follow. However, the techniques developed to prove Theorem~\ref{principal1} are flexible enough to handle perturbations of the atoms' masses as well. There are many possible variants of such perturbations; we are going to present one particular case below, and in Section \ref{generalization} we explain how each step of the proof has to be adapted to establish this more general version. Our aim is that this section serves as a guide to the reader who wants to explore other possible variants that can be addressed using the same methods. 

\medskip

Now, let us describe the variant that we have considered. First of all, for each $p\in X$, the atom with mass $\mu(p)$ can be decomposed into finitely many atoms. So, for each $p\in X$ we consider $n(p)\in \N$ (deterministic) and we consider the random perturbation 
\begin{equation}\label{general-perturbation}
\widetilde{\mu} :=\sum_{p\in X}\sum_{k=1}^{n(p)} \big(\alpha_{p,k}\,\mu(p)\big) \, \delta_{p + \xi_{p,k}},
\end{equation}
where we will assume the following hypotheses on the random vectors $\xi_{p,k}$ and the complex random variables $\alpha_{p,k}$: 

\begin{enumerate}[label=\textbf{(H\arabic*)},ref=H\arabic*]
    \item\label{eq:H1} \textbf{Independence for different values of $p$}: The families of random variables 
$$
\cP_p=\{\xi_{p,k},\alpha_{p,k}: 1\leq k \leq n(p)\}
$$ 
are independent for different values of $p$. However, for a fixed $p$, we allow dependencies among the variables with the same $p$, subject to hypothesis \textbf{(\ref{eq:H2})}. 

  \item\label{eq:H2} \textbf{Independence of each weight and its displacement}: For every $p\in X$ and $1\le k\le n(p)$, the random variable $\alpha_{p,k}$ and the random vector $\xi_{p,k}$ are independent.
  
  \item\label{eq:H3} \textbf{Tail control on the displacements}: There exists a random vector $\xi$ such that it controls the decay of the tails of the random vectors $\xi_{p,k}$ uniformly, more precisely, for every $\eta\geq 0$
\begin{equation}\label{Condition2 general}
\pe(\|\xi_{p,k}\|>\eta)\le\pe(\|\xi\|>\eta).
\end{equation}
We assume that the random vector $\xi$ satisfies the moment condition $\E(\|\xi\|^{d})<\infty.$

  \item\label{eq:H4} 
\textbf{The mass perturbation can not be arbitrary large}: The sum
\[
\sum_{k=1}^{n(p)} |\alpha_{p,k}| 
\]
is uniformly bounded on $p\in X$. 

  \item\label{eq:H5} \textbf{Asymptotic uniformity}: There exists a function $w:\R^d\to\C$ such that, for all $\la\in\R^d$,
\begin{equation}\label{Condition1 general}
\limsup_{\|p\| \to\infty}\left|  \e\left(\sum_{k=1}^{n(p)} \alpha_{p,k} \,e^{-2\pi i \langle \xi_{p,k}, \la \rangle}\right)-w(\la)\right|= 0.
\end{equation}

  \item\label{eq:H6}  \textbf{Finite variance}: In this general formulation, the variance hypothesis of Theorem \ref{principal1} is replaced by
\[
\sup_{p\in X} \ \e\left(\left(\sum_{k=1}^{n(p)} |\alpha_{p,k}|\lVert \xi_{p,k}\lVert\right)^2\right) <\infty. 
\]
\end{enumerate}

The reader may have the following particular case in mind: 

\begin{itemize}
\item The families $\cD_p=\{\xi_{p,k}: 1\leq k\leq n(p) \}$ are independent for different values of $p$, and all the random variables $\xi_{p,k}$ have the same distribution. We assume that the common distribution has finite variance and moment of order $d$.
\item The  families $\cC_p=\{\alpha_{p,k}: 1\leq k\leq n(p)\}$ are independent for different values of $p$ and with respect to all the displacements. Moreover, the random variables $\alpha_{p,k}$ take values in $[0,1]$ and satisfy 
$$
\sup_{p\in X}\ \sum_{k=1}^{n(p)} |\alpha_{p,k}| \leq 1 \qquad \mbox{and}\qquad  \lim_{\|p\|\to\infty} \e\left(\sum_{k=1}^{n(p)} \alpha_{p,k} \right)=1. 
$$
\end{itemize}

\begin{exa} 
An example of this situation in the plane is when, given $\eps>0$, for each atom $p\in X$, the mass $\mu(p)$ is spread uniformly on the $n(p)$-roots of the unity placed on the circle of radius $\eps$ and center $p$ and randomly shifted using a variable $\sigma_{p}$ uniformly distributed on $\{\|z\|=1\}$, so that for different values of $p$ these variables are independent. In this case
$$
\alpha_{p,k}=\frac{1}{n(p)}\qquad\mbox{and}\qquad \xi_{p,k}=\eps \ e^{2\pi i\frac{k}{n(p)}}\ \sigma_{p}.
$$ 
Note that all the variables $\xi_{p,k}$ are uniformly distributed on the circle $\{\|z\|=\eps\}$, they are independent for different values of $p$, but if we fix $p\in X$ the variables $\xi_{p,1},\ldots,\xi_{p,n(p)}$ are strongly correlated. 
\end{exa}

\medskip

Our second main theorem reads as follows. 

\medskip

\begin{teo}\label{principal2}

Let $\mu$ be a translation-bounded crystalline measure on $\R^d$, with $\mu$ and $\widehat{\mu}$ as in \textnormal{\textbf{%
\textcolor{black}{(}%
\hyperref[eq:N]{N}%
\textcolor{black}{)}}}. Let $\widetilde{\mu}$ be defined by
\eqref{general-perturbation}. Provided that hypotheses \textnormal{\textbf{(\ref{eq:H1})}-\textbf{(\ref{eq:H6})}}
hold, almost surely for every $\la\in\R^d$ we have that
\begin{equation}\label{mainresultgeneralized}
       \lim_{R\to\infty} \frac{1}{m_d(B_R)} \int_{B_R}e^{-2\pi i \langle x,\la \rangle} \,d\widetilde\mu(x) =  \left\{ \begin{array}{lcc}
             w(\lambda)A(\la) &   if  & \lambda\in X^* \\
             \\ 0 &  if  & \lambda\notin X^*
             \end{array}.
   \right.
    \end{equation}

\end{teo}

\subsubsection*{Diffraction}

The proof of Theorem \ref{principal1} requires, among other things,  analyzing the diffraction properties of the measure $\mu$, which is the other main goal of this paper. The obtained results extend those of Alon et al.~\cite[Section~2.3]{MR4841781} (see also \cite{MR4576026}), who obtain diffraction properties of crystalline measures in which all masses $\mu(p)$ are equal to one. In that setting, the support is a finite union of uniformly separated sets \cite[Lemma 2]{MR4633114}, in particular, the measure turns out to be translation-bounded. 

\medskip

Our main result in this direction is the following. The reader who is not familiar with diffraction theory is referred to Subsection \ref{Prelim difraccion}, where we recall the basic definitions and properties.

\begin{teor}\label{diffractionmu}
Let $\mu$ be a translation-bounded crystalline measure, with $\mu$ and $\widehat{\mu}$ as in \textnormal{\textbf{%
\textcolor{black}{(}%
\hyperref[eq:N]{N}%
\textcolor{black}{)}}}. The autocorrelation of \(\mu\) exists. Moreover, its diffraction is given by
\[
\widehat{\gamma_\mu}=\sum_{\lambda\in X^*}|A(\lambda)|^2\,\delta_\lambda \,.
\]
\end{teor}

Our proof is self-contained and relies only on elementary tools from Fourier analysis. The method simultaneously establishes the existence of the diffraction and determines its form. The same result can also be proved using methods from the theory of almost periodic functions (see Remark \ref{otroenfoque}). 
One immediate consequence is related to the hyperuniformity of $\mu$. Recall that, given a measure $\nu$ such that its autocorrelation $\gamma_\nu$ exists, it is \emph{hyperuniform} if the following condition holds
$$
\lim _{\eps \rightarrow 0} \frac{\widehat{\gamma}_\nu\left(B_\eps \setminus\{0\}\right)}{\eps^d}=0 .
$$
The reader is referred to \cite{MR4735950, MR3815253} for more information about hyperuniform distributions of points. A hyperuniform measure $\nu$ is called \emph{stealthy hyperuniform} if the quotient in the limit is identically zero for small enough values of $\eps$. Since $X^*$ is locally finite, as a direct consequence of the previous theorem is the following corollary.

\begin{cor}\label{es requeterecontra hiperuniforme}
If $\mu$ is a translation-bounded crystalline measure, 
then it is stealthy hyperuniform.
\end{cor}

\medskip

\subsection{Organization of the paper}

The paper is organized as follows. Section \ref{preliminares} is devoted to collecting some preliminaries on Radon measures and diffraction theory that will be used throughout the paper. In Section \ref{Difraccion of TB cristalline} we study the diffraction properties of translation-bounded crystalline measures. In Section  \ref{roadmap}, we describe the roadmap of the proof of Theorem \ref{principal1}. In particular, we explain how to split the problem into its structured, almost periodic, and random parts. Section \ref{random part} is devoted to studying the random part of the problem. It is divided into several subsections with the aim of simplifying its reading. With the same aim, in this section, we will focus on Theorem \ref{principal1}. This is done to keep the notation as simple as possible. One of the main results of this section is Theorem \ref{chuecageneralizada}, which is a kind of twisted law of large numbers. Finally, in Section \ref{generalization}, we then turn to its generalization, Theorem \ref{principal2}. In the appendix, we will prove a weighted law of large numbers, which is required in Sections \ref{random part} and \ref{generalization}.

\section{Preliminaries}\label{preliminares}

We keep the notation introduced in Subsection~\ref{framework}. The Schwartz class $\schw$ is endowed with the usual locally convex topology (see \cite[Chapter 7]{MR1157815}). Given $f\in \schw$, the \textit{Fourier transform} $\efe(f)=\widehat{f}$  is defined by
\begin{equation*}
        \widehat{f}(w) := \int_{\mathbb{R}^d} f(x)e^{-2\pi i \langle x,w \rangle} \, \deriv x \,,
    \end{equation*}
where $\langle \cdot,\cdot \rangle$  denotes the standard inner product of $\R^d$. The space of tempered distributions is denoted by  $\schwd$. Recall that, by duality, the \textit{distributional Fourier transform} of $T\in \schwd$ is defined as the unique tempered distribution $\widehat{T}$ such that
\begin{equation*}
    \widehat{T}(\fii) := T(\widehat{\fii}) \,, \quad \fii \in \schw \,.
\end{equation*}

In this paper we will only consider the Fourier transform in this distributional sense. The set of compactly supported continuous functions is denoted by $C_c(\mathbb{R}^d)$, and the set of compactly supported $C^\infty$ functions is denoted by $C_c^\infty(\R^d)$. By means of $m_d(B_R)$ we denote the Lebesgue measure of $B_R$. Finally, we will use the notation $f(x)\lesssim g(x)$ to express that there exists a constant $C$, independent of $f$ and $g$ such that $f(x)\leq C g(x)$ for every $x$ in the common domain of the functions $f$ and $g$.

\subsection{Radon measures}

 For a purely atomic Radon measure 
\begin{equation*}
  \nu = \sum_{p\in X} \nu(p)\delta_p,
\end{equation*}
the total variation is simply
\begin{equation*}
  |\nu| = \sum_{p\in X} |\nu(p)|\delta_p \,.
\end{equation*}
We say a (Radon) measure $\mu$ has polynomial growth if there exists $n\in\N$ such that
$$
|\mu|(B_R)\lesssim R^n.
$$
If $\mu$ has polynomial growth, then it defines a tempered distribution by integration. We refer to it in the same way, with a slight abuse of notation. For a \textit{translation-bounded}  (TB) measure $\mu$, we use the notation
\[
\|\mu\|_{TB}:=\sup_{x \in \mathbb{R}^d} |\mu|(B_1(x))  \,.
\]
This quantity defines a norm in the vector space of translation-bounded measures, and with this norm, it becomes a Banach space. Note that if $\mu$ is a TB measure, then there exists a purely dimensional constant $C_d$ such that
\begin{equation}\label{TBisPolynomial}
    |\mu|(B_R)\leq C_d\, \|\mu\|_{TB}\ R^d.
\end{equation}
Hence, any TB-measure has polynomial growth with degree equal to the dimension of the underlying space $\R^d$.

\medskip

Given a sequence $\{\mu_n\}_{n\geq 1}$ of measures, it converges \textit{weakly} to a measure $\mu$ if for every  $f\in C_c(\R^d)$
\begin{equation*}
    \lim_{n\to\infty} \mu_n(f) = \mu(f) .
\end{equation*}
If we have a family $\{\mu_R\}_{R>0}$, it converges \textit{weakly} to a measure $\mu$ if for any increasing sequence $\{R_n\}_{n\geq 1}$, the sequence $\{\mu_{R_n}\}_{n\geq 1}$ converges to $\mu$.

\medskip
On the other hand, the sequence $\{\mu_n\}_{n\geq 1}$ converges to $\mu$  \textit{in the distributional sense} if for every  $f\in \schw$
\begin{equation*}
    \lim_{n\to\infty} \mu_n(f) = \mu(f) .
\end{equation*}
A similar extension as before is applied to families of the form $\{\mu_R\}_{R>0}$. We will need the following lemma that relates weak convergence and convergence in distribution. For the sake of completeness, we also include a simple proof.
\begin{lem}\label{vagimpliesdist}
    Let $\{\nu_{t}\}$ be a family of Radon measures that have polynomial growth \textit{uniformly}, more precisely, there exists $C>0$ and $n\in\N$ such that
    \begin{equation*}
        |\nu_{t}|(B_{R})\leq C\, {R}^n \,,
    \end{equation*}
    for all $t$ and $R\geq 1$. Assume that this family weakly converges to a measure $\nu$ which also has polynomial growth
    \begin{equation*}
        |\nu|(B_{R})\lesssim {R}^m \,.
    \end{equation*}
    Then, the family  $\{\nu_{t}\}$ also converges to $\nu$ in the distributional sense. 
\end{lem}
\begin{proof}
    Since the measures $\{\nu_{t}\}$ and $\nu$ have polynomial growth, they define tempered distributions by integration. Moreover, there exists some constants such that 
    \begin{align*}
        |\nu_{t}(f)|\leq C'\lVert f\lVert_\a \  (\mbox{for every $t $}) \qquad \mbox{and}\qquad 
        |\nu(f)|\leq C''\lVert f\lVert_\beta \,,
    \end{align*}
    for each $f\in\schw$, and some norms $\lVert\cdot\lVert_\a$, $\lVert\cdot\lVert_\beta$ on $\schw$. Given $\fii\in \schw$ and $\eps>0$, let $\psi\in C_c^\infty(\R^d)$ such that
    \begin{equation*}
        \lVert \fii-\psi \lVert_\a \leq \frac{\eps}{3C'} \,\quad \text{and}\,\quad \lVert \fii-\psi \lVert_\beta \leq \frac{\eps}{3C''} \,.
    \end{equation*}
    Let $t_0$ such that 
    \begin{equation*}
        |\nu_{t}(\psi)-\nu(\psi)|\leq \frac{\eps}{3} \,.
    \end{equation*}
    for all $t \geq t_0$. Then if $t\geq t_0$:
    \begin{equation*}
        |\nu_{t}(\fii)-\nu(\fii)|\leq |\nu_{t}(\fii-\psi)| + |\nu_{t }(\psi)-\nu(\psi)| + |\nu(\psi-\fii)| \leq \eps \,.
    \end{equation*}
\end{proof}

\medskip

The following lemma will also be useful in the sequel. The proof is based on standard covering arguments. 

\begin{lem}\label{mugrowth}
Let $\nu$ be a translation-bounded measure on $\R^d$. Then there exists a constant \(C=C(d)>0\) such that
\begin{equation*}
     |\nu|(B_{R_2}(q))-|\nu|(B_{R_1}(q)) \leq C \|\nu\|_{TB}\ R_2^{d-1} \max\{1,R_2-R_1\}\,, \quad q\in \R^d \,,
\end{equation*}
whenever $R_2>R_1\geq 0$.
\end{lem}


\subsection{Diffraction theory}\label{Prelim difraccion}

Recall that given two finite measures $\nu_1$ and $\nu_2$, their convolution is the measure $\nu_1\ast \nu_2$ that acts on a function $f\in C_c(\R^d)$ as
$$
(\nu_1 \ast \nu_2)(f) := \int_{\R^d \times \R^d} f(x+y) \,d\nu_1(x)d\nu_2(y) \,.
$$
This can be extended to characteristic functions of Borel measurable sets by standard arguments. The convolution is also well-defined when one of the measures $\nu_j$ is translation-bounded. 

\medskip

Given a function $f\in C_c(\R^d)$ 
$$
\widetilde{f}(x):=\overline{f(-x)}.
$$ 
This can be extended to measures by defining $\widetilde{\mu}$ on a function $f\in C_c(\R^d)$ as follows
$$
\widetilde{\mu}(f):=\overline{\mu(\widetilde{f})} \,.
$$

\medskip

Now, recall that the autocorrelation of $f\in L^1(\R^d)$ is defined by
$$
\gamma_f= f*\widetilde{f}.
$$
This can be extended to a finite measure. However, since the convolution of two arbitrary measures is not defined, the autocorrelation of a Radon measure is not direct, and it may not exist. Consider a measure $\mu$. Given $R>0$, define the finite measure 
$$
\mu_R(A):=\mu(A\cap B_R) \,.
$$
The \textit{volume-averaged convolution} is defined as
$$
\gamma_\mu^{(R)} := \frac{\mu_R \ast \widetilde{\mu_R}}{m_d(B_R)} \,.
$$

\medskip

\begin{exa}
Let $\mu$ be a purely atomic measure, that is
\begin{equation*}
  \mu = \sum_{p\in X} \mu(p)\delta_p,
\end{equation*}
and assume that its support $X$ is locally finite. Then, 
\begin{equation*}
  \mu_R = \sum_{p\in X_R} \mu(p)\delta_p \qquad\mbox{and}\qquad \widetilde{\mu_R}=\displaystyle\sum_{p\in X_R} \overline{\mu(p)}\delta_{-p} \,.
\end{equation*}
Therefore,  we get that
   \begin{equation*}
    \gamma_{\mu}^{(R)} = \frac{1}{m_d(B_R)} \sum_{p,q\in X_R} \mu(p)\overline{\mu(q)} \delta_{p-q} \,.
\end{equation*}
\end{exa}

\medskip

If $\mu$ is translation-bounded, then the family $\{\mu_{R}\}_{R>0}$ is precompact in the weak topology. Moreover, they are uniformly bounded with respect to the norm $\|\cdot\|_{TB}$. Hence, there exists $C>0$ such that for every $R>0$ and $r\geq 1$
   \begin{equation}\label{eq TB uniformly d growth}
        |\mu_{R}|(B_{r})\leq C\, {r}^d.
    \end{equation}
These properties are inherited by the family $\{\gamma_\mu^{(R)} \}$. In particular, by the properties of the weak topology and Lemma \ref{vagimpliesdist}, this family has sequential limits with respect to the weak topology, and with respect to the convergence in distribution as well (see \cite[Proposition 9.1]{MR3136260}). 

\begin{fed}[Baake-Grimm \cite{MR3136260}, Hof \cite{MR1328260}]
Given a Radon measure $\mu$, if the limit 
$$
\lim_{R\to \infty} \gamma_\mu^{(R)}
$$
exists in the weak sense, it is called the autocorrelation of the measure $\mu$, and it is denoted by $\gamma_\mu$.

\end{fed}

If a measure $\mu$ is TB, the autocorrelation exists if the family $\{\gamma_\mu^{(R)} \}$ has a unique sequential limit in the weak sense. By the properties of the family $\{\gamma_\mu^{(R)} \}$ when $\mu$ is translation-bounded, this unique limit is also the unique limit in the distributional sense. In any case, it is not difficult to prove that any weak limit $\la$ of the family $\{\gamma_\mu^{(R)} \}$ is completely positive, in the sense that $\la(f\ast \widetilde{f})\geq 0$ for every $f\in \schw$. Therefore, by  Bochner's theorem, the Fourier transform of each $\la$ (in the distributional sense) is a positive measure. 
    
 \begin{fed}[Baake-Grimm \cite{MR3136260}, Hof \cite{MR1328260}]
If the autocorrelation of a TB measure $\mu$ exists, then the \textit{diffraction measure} associated with $\mu$ is the Fourier transform of $\gamma_\mu$ (in the distributional sense). 
\end{fed}   

 The interested reader is referred to \cite{MR3136260} and \cite{MR1328260} for the proofs of all the aforementioned results and more information. 

\medskip

\section{Diffraction of translation-bounded crystalline measures}\label{Difraccion of TB cristalline}

In this section, we will study the diffraction properties of translation-bounded crystalline measures. We will start proving some consequences of the generalized Poisson formula associated with a crystalline measure, which will also be useful in the sequel. The second part of the section is devoted to proving Theorem \ref{diffractionmu}.

\subsection{First consequences of the generalized Poisson summation formula}

The formula $\eqref{poisson}$ can be rewritten in the following way.

\begin{prop}\label{notPoisson}
Let $\mu$ be a translation-bounded crystalline measure, with 
\[
\mu=\sum_{p\in X}\mu(p)\delta_p,
\qquad
\widehat{\mu}=\sum_{\lambda\in X^*}A(\lambda)\delta_\lambda \,.
\]
Then, for every $\la\in\R^d$
\begin{equation*}
   \lim_{R\to\infty} \frac{1}{m_d(B_R)} \sum_{p\in X_R} e^{-2\pi i \langle p,\lambda \rangle}\mu(p) = \left\{ \begin{array}{lcc}
             A(\lambda) &   if  & \lambda\in X^* \\
             \\ 0 &  if  & \lambda\notin X^*
             \end{array}
   \right. \,.
\end{equation*}
\end{prop}

\bigskip

\noindent The proof of this proposition is based on the following two lemmas. 

\begin{lem}
Let \(\varphi\in\mathcal{S}(\mathbb{R}^d)\), and define $\fii_R(x):=\fii(Rx)$ for every $R>0$. For any \(T\in\mathcal{S}'(\mathbb{R}^d)\) satisfying \(x_0\notin\operatorname{supp}(T)\), we have
\[
\lim_{R\to\infty} T(\varphi_R(\cdot- x_0))=0 \,.
\]
\end{lem}

\begin{proof}
For simplicity, let us assume that $x_0=0$. The general case follows essentially the same arguments. Choose \(\delta>0\) such that
$B_{2\delta}\cap\operatorname{supp}(T)=\emptyset$, and take
$\rho\in C_c^\infty(\R^d)$ with $\rho\equiv1$ on $B_\delta$ and $\operatorname{supp}(\rho)\subseteq B_{2\delta}$.
By Leibniz's rule and the rapid decay of $\varphi$, for every
pair of multi-indices $\alpha,\beta$ and every
$N>|\alpha|+|\beta|$, we have, for $R\geq 1$
\[
\sup_{x\in\mathbb R^d}
\left|x^\alpha\partial^\beta
\bigl((1-\rho(x))\varphi_R(x)\bigr)\right|
\le C_{\alpha,\beta,N}R^{|\beta|-N}
\xrightarrow[R\to\infty]{} 0.
\]
Thus $(1-\rho)\varphi_R\to0$ in the Schwartz topology,
that is, in every Schwartz seminorm.
Since $\operatorname{supp}(\rho\varphi_R)\subseteq B_{2\delta}$,
we have $T(\rho\varphi_R)=0$. By continuity of $T$,
\[
T(\varphi_R)=T\bigl((1-\rho)\varphi_R\bigr)\xrightarrow[R\to\infty]{} 0 \,.
\] 
\end{proof}

\begin{lem}\label{limfii}
Let \(\varphi\in\mathcal{S}(\mathbb{R}^d)\) and \(t,\gamma\in\mathbb{R}^d\). Then,
\begin{equation*}
    \lim_{R\to\infty} \frac{1}{R^d} \sum_{p\in X} \widehat{\fii}\left(\frac{p-t}{R}\right)e^{-2\pi i \langle p,\gamma \rangle}\mu(p) = \left\{ \begin{array}{lcc}
             A(\gamma)\fii(0) &   if  & \gamma\in X^* \\
             \\ 0 &  if  & \gamma\notin X^*
             \end{array}
   \right. \,.
\end{equation*}
\end{lem}

\begin{proof} 
Let \(T\in\mathcal{S}'(\mathbb{R}^d)\) be given by
\[
T(f):=\sum_{\lambda\in X^*\setminus\{\gamma\}}A(\lambda)e^{2\pi i\langle t,\lambda\rangle}f(\lambda) \,,\quad f\in\mathcal{S}(\mathbb{R}^d) \,,
\]
and let us extend the function $A$, a priori only defined on $X^*$, as zero outside the set $X^*$. Using the Poisson summation formula and $\fii_R(x)=\fii(Rx)$ as before we can write
\begin{align*}
\lim_{R\to\infty}\frac{1}{R^d}\sum_{p\in X}\widehat{\varphi}\Bigl(\frac{p-t}{R}\Bigr)e^{-2\pi i\langle p,\gamma\rangle}\mu(p)
&= e^{-2\pi i\langle \gamma,t\rangle}\lim_{R\to\infty}\sum_{\lambda\in X^*}A(\lambda)\,\varphi(R(\lambda-\gamma))e^{2\pi i\langle \lambda,t\rangle}\\[1mm]
&= e^{-2\pi i\langle \gamma,t\rangle}\Bigl(\lim_{R\to\infty}T(\varphi_R(\cdot-\gamma))+A(\gamma)\varphi(0)e^{2\pi i\langle \gamma,t\rangle}\Bigr) \,.
\end{align*}
By the previous lemma the first term vanishes, yielding the desired result.
\end{proof}

\medskip

\begin{proof}[Proof of Proposition~\ref{notPoisson}]
For a given \(0<\delta<1\), choose \(\varphi_\delta\in\mathcal{S}(\mathbb{R}^d)\) such that its Fourier transform \(\widehat{\varphi_\delta}\) is real, radial, decreasing, and satisfies
\[
\widehat{\varphi_\delta}\equiv1\quad\text{on } B_1\,,\qquad \widehat{\varphi_\delta}\equiv0\quad\text{on } B_{1+\delta}^c \,.
\]
Then,
\[
\frac{1}{R^d}\sum_{p\in X_R}e^{-2\pi i\langle p,\lambda\rangle}\mu(p)
=\frac{1}{R^d}\sum_{p\in X}e^{-2\pi i\langle p,\lambda\rangle}\mu(p)\,\widehat{\varphi_\delta}\Bigl(\frac{p}{R}\Bigr)
-\frac{1}{R^d}\sum_{p\in X_{R+R\delta}\setminus X_R}e^{-2\pi i\langle p,\lambda\rangle}\mu(p)\,\widehat{\varphi_\delta}\Bigl(\frac{p}{R}\Bigr) \,.
\]
On the one hand, for \(R\ge1/\delta\) (so that Lemma \ref{mugrowth} applies), we have
\begin{align*}
\left|\frac{1}{R^d}\sum_{p\in X_{R+R\delta}\setminus X_R}e^{-2\pi i\langle p,\lambda\rangle}\mu(p)\,\widehat{\varphi_\delta}\Bigl(\frac{p}{R}\Bigr)\right|
&\le\frac{|\mu|(X_{R+R\delta})-|\mu|(X_R)}{R^d} \\
&\lesssim\frac{(R+R\delta)^{d-1}(R\delta)}{R^d} = (1+\delta)^{d-1}\delta \,.
\end{align*}
On the other hand, 
\[
\varphi_\delta(0)=\int_{\mathbb{R}^d}\widehat{\varphi_\delta}(w)\,dw
=m_d(B_1)+\int_{B_{1+\delta}\setminus B_1}\widehat{\varphi_\delta}(w)\,dw \,.
\]
Hence, \(\varphi_\delta(0)\to m_d(B_1)\) as \(\delta\to0\). Therefore, by Lemma~\ref{limfii}, the proof is complete.
\end{proof}

\medskip

\begin{cor}\label{A lambda unif acotados}
Let $\mu$ be a translation-bounded crystalline measure, with \[
\mu=\sum_{p\in X}\mu(p)\delta_p,
\qquad
\widehat{\mu}=\sum_{\lambda\in X^*}A(\lambda)\delta_\lambda \,.
\]
Then, 
$$
\sup_{\la\in X^*} |A(\la)|\leq C_d\|\mu\|_{TB},
$$
where $C_d$ is a purely dimensional constant.
\end{cor}

\smallskip

Indeed, by Proposition \ref{notPoisson}
\begin{equation}\label{cotaAlambda}
|A(\lambda)|=\lim_{R\to\infty} \left|\frac{1}{m_d(B_R)} \sum_{p\in X_R} e^{-2\pi i \langle p,\lambda \rangle}\mu(p)\right|\le\liminf_{R\to\infty}\frac{|\mu|(B_R)}{m_d(B_R)}\,,
\end{equation}
and by $\eqref{TBisPolynomial}$, the right-hand side can be bounded by a purely dimensional constant times $\|\mu\|_{TB}$. Another consequence of \eqref{cotaAlambda} is the following lower bound for the growth of the total variation of $\mu$.

\begin{cor}\label{muinfinita}
Let $\mu$ be a (not identically zero) translation-bounded crystalline measure. For large enough $R>0$:
\[
R^d \lesssim  |\mu|(B_R) \,.
\]
\end{cor}

\begin{proof}
    Since \(\mu\neq0\), its Fourier transform is nonzero, so there exists \(\lambda_0\in X^*\) such that \(A(\lambda_0)\neq0\). By $\eqref{cotaAlambda}$,
$$
\liminf_{R\to\infty}\frac{|\mu|(B_R)}{m_d(B_R)}
\geq |A(\lambda_0)|>0\,,
$$
which proves the claim.

\end{proof}

\subsection{Formula for the diffraction of the crystalline measure}

This subsection is devoted to prove Theorem \ref{diffractionmu}. More precisely, we have to see that for all $f\in C_c(\R^d)$

\begin{equation}\label{existsdiffrac}
    \lim_{R\to\infty} \frac{1}{m_d(B_R)} \sum_{p,q\in X_R} \mu(p)\overline{\mu(q)} f(p-q)
\end{equation}
exists. As a byproduct of our proof, we will get that for any $f\in \schw$
\begin{equation}\label{eqdiffractionmu}
    \lim_{R\to\infty} \frac{1}{m_d(B_R)} \sum_{p,q\in X_R} \mu(p)\overline{\mu(q)} \widehat{f}(p-q) = \sum_{\la\in X^*} |A(\la)|^2 f(\la) \,.
\end{equation}

The first step is the following lemma, which is a particular case of \cite[Lemma 1.2]{MR1798991}. 
\begin{lem}
Given two purely atomic translation-bounded measures
\begin{equation*}
        \nu_1 = \sum_{p\in X} \nu_1(p)\delta_p \,, \quad\mbox{and}\quad \nu_2 = \sum_{p\in X} \nu_2(p)\delta_p,
\end{equation*} 
for every $\varphi\in C_c(\R^d)$ 
    \begin{equation}\label{martin}
        \lim_{R\to\infty} \frac{1}{m_d(B_R)}\sum_{q\in X_R}\sum_{p\in X\setminus X_R}\nu_1(p)\overline{\nu_2(q)}\varphi(p-q)
=0 \,.
   \end{equation}
\end{lem}

Since our definition of Fourier transform is in the distributional sense, we have to extend this result to functions in $\schw$.

\begin{prop}\label{sevuelveausar} 
  The limit $\eqref{martin}$ is valid for $\fii\in\schw$ too.
\end{prop}

\begin{proof}
Let's take
\begin{equation*}
    \eta_R := \frac{1}{m_d(B_R)}\sum_{q\in X_R} \sum_{p\in X\setminus X_R} \nu_1(p)\overline{\nu_2(q)} \delta_{p-q} \,, \quad R\geq 1 \,.
\end{equation*}
The family $\{\eta_R\}$ is \textit{uniformly} translation-bounded, in the sense that there exists a $C>0$ so that
\begin{equation*}
    \sup_{R\geq 1}\|\eta_R\|_{TB} \leq \tilde{C} \,.
\end{equation*}
Indeed, given any $q_0\in\R^d$
\begin{align*}
\sum_{q\in X_R}\sum_{p\in X\setminus X_R}|\nu_1(p)||\nu_2(q)||\mathbbm{1}_{B_1(q_0)}(p-q)|
&=\sum_{q\in X_R}|\nu_2(q)|\sum_{p\in (X\setminus X_R)\cap B_{1}(q+q_0)}|\nu_1(p)| \\[1mm]
&\leq \sum_{q\in X_R} |\nu_2(q)|  |\nu_1|(B_1(q+q_0)) \lesssim R^{d} \,,
\end{align*}

where we have used $\eqref{TBisPolynomial}$. As a consequence, the family $\{\eta_R\}$ has \text{uniformly} polynomial growth. With $\eqref{martin}$ and Lemma \ref{vagimpliesdist}, we obtain the result.
\end{proof}

\medskip

Before starting the proof of Theorem \ref{diffractionmu}, note that for any $f\in \schw$
\begin{equation*}
    \lim_{R\to\infty} \frac{1}{m_d(B_R)} \sum_{q\in X_R}\sum_{p\in X_R} \mu(p)\overline{\mu(q)} \widehat{f}(p-q) =     \lim_{R\to\infty} \frac{1}{m_d(B_R)} \sum_{q\in X_R}\sum_{p \in X} \mu(p)\overline{\mu(q)} \widehat{f}(p-q) \,.
\end{equation*}
So, as a direct consequence of the Poisson summation formula and Proposition \ref{notPoisson}, we have that for any $f\in C_c^\infty(\mathbb{R}^d)$:
    \begin{align*}
         \lim_{R\to \infty} \frac{1}{m_d(B_R)}\sum_{q\in X_R}\sum_{p\in X} \mu(p) \overline{\mu(q)}\widehat{f}(p - q) &=\lim_{R\to \infty} \frac{1}{m_d(B_R)}\sum_{q\in X_R} \overline{\mu(q)} \sum_{\la\in X^*} A(\la)e^{2\pi i \langle \la,q \rangle}f(\la)  \\[1mm]
         &=\sum_{\la\in X^*} A(\la) f(\la) \lim_{R\to \infty} \frac{1}{m_d(B_R)} \sum_{q\in X_R} e^{2\pi i \langle \la, q \rangle} \overline{\mu(q)}  \\[1mm]
         &=\sum_{\lambda\in X^*} |A(\lambda)|^2f(\lambda)  \,.
\end{align*} 

We are now ready to go into the proof of Theorem \ref{diffractionmu}. 

\begin{proof}[Proof of Theorem \ref{diffractionmu}]
As shown in \cite[Proposition 9.1]{MR3136260} the family of measures 
\begin{equation*}
    \gamma_{\mu}^{(R)} = \frac{1}{m_d(B_R)} \sum_{p,q\in X_R} \mu(p)\overline{\mu(q)} \delta_{p-q} \,,\quad R\geq 1
\end{equation*}
is \textit{uniformly} translation-bounded, that is
\begin{equation*}
   \sup_{R>0} \|\gamma_{\mu}^{(R)}\|_{TB}\leq C \,,
\end{equation*}
and it is precompact in the weak topology. So, given a sequence $\{\gamma_{\mu}^{(R_k)}\}_{k\geq 1}$, there is a subsequence $\{\gamma_{\mu}^{(R_{k_\ell})}\}_{\ell\geq 1}$ and a measure $\gamma$ such that
\begin{equation*}
        \gamma_{\mu}^{(R_{k_\ell})} \xrightarrow[\ell\to\infty]{\text{weakly}} \gamma.
    \end{equation*} 
Moreover, it is not difficult to see that $\|\gamma\|_{TB}\leq C$ as well. 
Note that, by Lemma \ref{vagimpliesdist}, the convergence to $\gamma$ not only is in the weak topology, but also in the distributional sense. This implies that
\begin{equation*}
        \efe\bigl(\gamma_{\mu}^{(R_{k_\ell})}\bigr) \xrightarrow[\ell\to\infty ]{S'} \efe(\gamma) \,.
\end{equation*} 
Since $\{\efe\bigl(\gamma_{\mu}^{(R_{k_\ell})}\bigr) \}$ are positive measures,  $\efe(\gamma)$ is a positive measure too. Then, by \cite[Proposition 3.3]{MR1328260}, $\efe(\gamma)$ is translation-bounded. We have shown that for \(f\in C_c^\infty(\mathbb{R}^d)\) 
\begin{equation}\label{step4}
\lim_{R\to\infty}\frac{1}{m_d(B_R)}\sum_{q\in X_R}\sum_{p\in X_R}\mu(p)\overline{\mu(q)}\widehat{f}(p-q)
=\sum_{\lambda\in X^*}|A(\lambda)|^2\,f(\lambda) \,.
\end{equation}
Hence, for any $f\in C_c^\infty(\R^d)$ we have
    \begin{equation}\label{igualdadcc}
        \efe(\gamma)(f)= \sum_{\la\in X^*} |A(\la)|^2 f(\la) \,. 
    \end{equation}
This implies that the measure
    \begin{equation*}
        \sum_{\la\in X^*} |A(\la)|^2 \delta_{\la}
    \end{equation*}
is translation-bounded, and therefore it defines a tempered distribution. Hence, the equality $\eqref{igualdadcc}$ extends to any function in $\schw$. By the uniqueness of the Fourier transform, $\gamma$ does not depend on the original sequence. So the autocorrelation of \(\mu\) exists and is equal to \(\gamma\). Moreover, its diffraction is the TB-measure
    \begin{equation*}
    \widehat{\gamma_{\mu}} = \sum_{\la\in X^*} |A(\la)|^2 \delta_{\la} \,,
    \end{equation*}
which completes the proof.

\end{proof}

Set 
\[
\mu_\lambda:=\sum_{p\in X} e^{-2\pi i\langle \lambda,p\rangle}\mu(p)\,\delta_p \,, \qquad \la\in\R^d \,.
\]
The following corollary will be useful in Section \ref{random part}.

\begin{cor}\label{diffracdeterm}
For any $\la\in\R^d$, the autocorrelation of \(\mu_\lambda\) exists and its diffraction is given by
\[
\widehat{\gamma_{\mu_\lambda}}=\sum_{u\in X^*-\lambda}|A(u+\lambda)|^2\,\delta_u \,.
\]
\end{cor}

\begin{obs}
    In \cite{boyvalenkov2025growthmassescrystallinemeasures}, the authors show that $\sum_{\la\in X^*} |A(\la)|^2 \delta_{\la}$ is translation-bounded using a different approach.
\end{obs}

\begin{obs}\label{otroenfoque}
By \cite[Lemma 3]{MR4633114}, for each $\la\in\R^d$ the translation-bounded crystalline measure $\mu_\la$ is almost periodic. The reader is referred to \cite{MR4633114, MR2876415, MR3791850} and the references therein for more details on almost periodic functions, measures and distributions. Using this theory, alternative proofs of Proposition \ref{notPoisson} and Theorem \ref{diffractionmu} can be obtained from the results developed in \cite{MR3939563}.

\end{obs}

\section{Splitting into structured, almost periodic and random parts}\label{roadmap}

The aim of this section is to prove Theorem \ref{principal1}. Recall that $\mu$ is a translation-bounded crystalline measure, with 
\[
\mu=\sum_{p\in X}\mu(p)\delta_p,
\qquad
\widehat{\mu}=\sum_{\lambda\in X^*}A(\lambda)\delta_\lambda \,.
\]
As in Theorem \ref{principal1}, we consider the random perturbation of $\mu$, which we denote $\widetilde\mu$, and it is defined by
$$
\widetilde{\mu} :=\sum_{p\in X} \mu(p) \delta_{p + \xi_p}.
$$
In view of the remark following Theorem~\ref{principal1}, throughout
the proof we work under the more general assumptions {\hypersetup{linkcolor=azul}%
\textnormal{\textbf{(\ref{eq:A1})}}}-{\hypersetup{linkcolor=azul}%
\textnormal{\textbf{(\ref{eq:A2})}}}. Thus $ \{\xi_{p}\}_{p\in X}$ is an independent family of random vectors taking values in $\R^d$ such that there exists a function $w:\R^d\to\C$ such that
\[
\limsup_{\|p\|\to\infty} \left| \E(e^{-2\pi i \langle \xi_p ,\la \rangle })-w(\la) \right|=0,
\]
and there exists a random vector $\xi$ such that it controls uniformly the tail's decay of each random vector $\xi_p$. More precisely, for every $\eta\geq 0$ 
\[
\pe(\|\xi_p\|>\eta)\le\pe(\|\xi\|>\eta).
\]
We assume that the random vector $\xi$ also satisfies the moment condition 
$$
\E(\|\xi\|^{\max\{2,d\}})<\infty.
$$
For this section, we only need the $d$-moment to be finite. The second moment and the independence are not necessary until the next section. 

\medskip

Notice that
\begin{equation*}
    \frac{1}{m_d(B_R)} \int_{B_R}e^{-2\pi i \langle x,\la \rangle} \,d\widetilde\mu(x) = \frac{1}{m_d(B_R)}\sum_{p\in X:p+\xi_p\in B_R} e^{-2\pi i\langle p+\xi_p,\lambda\rangle}\mu(p) \,.
\end{equation*}

\medskip

The first step to prove Theorem \ref{principal1} is to divide the problem into three parts: one structured, one almost periodic, and the other random. The next result is crucial for making this division. Recall that given a set $Y\subseteq \R^d$, by means of $Y_R$ we denote the intersection $Y\cap B_R$.

\begin{pro}\label{key para dividir}
Almost surely,
\[
\sup_{\lambda\in\mathbb{R}^d}\Biggl|\frac{1}{m_d(B_R)}\sum_{p\in X:p+\xi_p\in B_R} e^{-2\pi i\langle p+\xi_p,\lambda\rangle}\mu(p)-\frac{1}{m_d(B_R)}\sum_{p\in X_R} e^{-2\pi i\langle p+\xi_p,\lambda\rangle}\mu(p)\Biggr|
\xrightarrow[R\to\infty]{} 0 \,.
\]
\end{pro}

We will prove this proposition in the next subsection. Now, let us see how this proposition permits us to divide the problem into three parts. First of all, as a consequence of this proposition, it is enough to prove
$$
       \lim_{R\to\infty} \frac{1}{m_d(B_R)} \sum_{p\in X_R} e^{-2\pi i\langle p+\xi_p,\lambda\rangle}\mu(p) =  \left\{ \begin{array}{lcc}
             w(\la)A(\lambda) &   if  & \lambda\in X^* \\
             \\ 0 &  if  & \lambda\notin X^*
             \end{array}.
   \right.
$$

A simple algebraic manipulation leads to the following splitting of the left-hand side

\begin{align*}
\frac{1}{m_d(B_R)}\sum_{p\in X_R}e^{-2\pi i\langle p+\xi_p,\lambda\rangle}\mu(p)
&=(\mbox{Structured part})+(\mbox{Almost periodic part})+(\mbox{Random part}),
\end{align*}
where 
\begin{align*}
(\mbox{Structured part})&=w(\la)\frac{1}{m_d(B_R)}\sum_{p\in X_R}e^{-2\pi i\langle p,\lambda\rangle}\mu(p),\\
(\mbox{Almost periodic part})&=\frac{1}{m_d(B_R)}\sum_{p\in X_R}e^{-2\pi i\langle p,\lambda\rangle}\mu(p) \, \big(\E(e^{-2\pi i \langle \xi_p ,\la \rangle })-w(\la) \big),\\
 (\mbox{Random part})  &=\frac{1}{m_d(B_R)}\sum_{p\in X_R}\Bigl[e^{-2\pi i\langle p,\lambda\rangle}\mu(p)\Bigl(e^{-2\pi i\langle \xi_p,\lambda\rangle}-\mathbb{E}\left(e^{-2\pi i\langle \xi_p,\lambda\rangle}\right)\Bigr)\Bigr] \,.
\end{align*}

Note that, by Proposition \ref{notPoisson}
\begin{equation*}
    \lim_{R\to\infty} \frac{1}{m_d(B_R)} \sum_{p\in X_R} e^{-2\pi i \langle p,\lambda \rangle}\mu(p) = \left\{ \begin{array}{lcc}
             A(\lambda) &   if  & \lambda\in X^* \\
             \\ 0 &  if  & \lambda\notin X^*
             \end{array}
   \right. \,.
\end{equation*}
Therefore, the limit of the structured part is precisely the limit we are seeking. On the other hand,
\begin{align*}
\left|\frac{1}{m_d(B_R)}\sum_{p\in X_R}e^{-2\pi i\langle p,\lambda\rangle}\mu(p) \, \big(\E(e^{-2\pi i \langle \xi_p ,\la \rangle })-w(\la) \big)\right|&\leq \frac{1}{m_d(B_R)}\sum_{p\in X_R}|\mu(p)| \, |\E(e^{-2\pi i \langle \xi_p ,\la \rangle })-w(\la) |\\
&\lesssim\frac{1}{|\mu|(B_R)} \sum_{p\in X_R}|\mu(p)| \, |\E(e^{-2\pi i \langle \xi_p ,\la \rangle })-w(\la) |,
\end{align*}
where for the last inequality we have used that $\mu$ is translation-bounded. Now, since the random vectors $\xi_p$ satisfy condition {\hypersetup{linkcolor=azul}%
\textnormal{\textbf{(\ref{eq:A1})}}}, 
\begin{align*}
\lim_{R\to\infty}\frac{1}{|\mu|(B_R)} \sum_{p\in X_R}|\mu(p)| \, |\E(e^{-2\pi i \langle \xi_p ,\la \rangle })-w(\la) |=0 \,.
\end{align*}
Therefore, to finish the proof of Theorem \ref{principal1}, we have to prove that the limit of the random part is zero. The next section is devoted to this problem. Before that, we will prove Proposition \ref{key para dividir}. 

\subsection{Proof of Proposition \ref{key para dividir}}

We have to prove that, almost surely
\[
\sup_{\lambda\in\mathbb{R}^d}\Biggl|\frac{1}{m_d(B_R)}\sum_{p\in X:p+\xi_p\in B_R} e^{-2\pi i\langle p+\xi_p,\lambda\rangle}\mu(p)-\frac{1}{m_d(B_R)}\sum_{p\in X_R} e^{-2\pi i\langle p+\xi_p,\lambda\rangle}\mu(p)\Biggr|
\xrightarrow[R\to\infty]{} 0 \,.
\]

To prove this, we have to estimate the difference between both quantities due to those points of $X$ inside $B_R$ that go out of $B_R$ when they are perturbed, and those points of $X$ that are outside $B_R$ and go inside the ball when they are perturbed. As the number of points of $X_R$ can grow faster than $R^d$, the proof has to involve the total variation of $\mu$ instead of the number of points of $X$.

\medskip

\begin{lem}\label{adentroafuera}
Almost surely,
\[
\lim_{R\to\infty}\frac{1}{R^d}|\mu|(\{p\in X:\|p\|\le R, \|p+\xi_p\|>R\})=0 
\]
and
\[
\lim_{R\to\infty}\frac{1}{R^d}|\mu|(\{p\in X:\|p\|>R, \|p+\xi_p\|\le R\})=0 \,.
\]
\end{lem}

\begin{proof}
Notice that
\[
A_R:= \{p\in X:\|p\|\le R, \|p+\xi_p\|>R\} \subseteq  \{p\in X: R-\|\xi_p\|<\|p\|\le R\}\,. 
\]
Hence
\[
|\mu|(A_R)\leq \sum_{p\in X_R} |\mu(p)|\mathbbm{1}_{\{\lVert p\rVert > R - \lVert \xi_p\rVert\}} \,.
\]
So
\[
\e\left(|\mu|(A_R)\right) \leq \sum_{p\in X_{R}} |\mu(p)| \e(\mathbbm{1}_{\{\lVert p\rVert > R - \lVert \xi\rVert\}})  = \mbox{(I)} + \mbox{(II)} \,,
\]
where
\begin{align*}
\mbox{(I)} &= \sum_{p\in X_{R}} |\mu(p)| \left[\e(\mathbbm{1}_{\{\lVert p\rVert > R - \lVert \xi\rVert\}} \,\,|\,\, \|\xi\|\leq R ) \,\pe(\|\xi\|\leq R )\right]  \\
&=\pe(\|\xi\|\leq R )\,\e\left(\left.\sum_{p\in X_{R}\setminus X_{R-\|\xi\|}} |\mu(p)| \right| \|\xi\|\leq R\right)\\[2mm] 
&\lesssim R^{d-1}\e(\|\xi\| + 1) 
\end{align*}
and
\begin{align*}
\mbox{(II)} &= \sum_{p\in X_{R}} |\mu(p)| \left[\e(\mathbbm{1}_{\{\lVert p\rVert > R - \lVert \xi\rVert\}} \,\,|\,\,\|\xi\|> R ) \pe(\|\xi\|> R )\right]  \\
&=\pe(\|\xi\|> R )|\mu|(X_{R}) \lesssim R^d \ \frac{\e(\|\xi\|)}{R} \,,
\end{align*}
where we have used Markov's inequality.
We get
\[
\e\left(\frac{1}{R^d} |\mu|(A_R)\right) \lesssim \frac{1}{R}\left(\e(\lVert \xi\Vert) + 1\right) \,.
\]
Since 
\[
\e\left(\sum_{n\ge1}\frac{1}{(n^2)^d} |\mu|(A_{n^2})\right) = \sum_{n\ge1}\e\left(\frac{1}{(n^2)^d} |\mu|(A_{n^2})\right)<\infty \,,
\]
then 
\[
\sum_{n\ge1}\frac{1}{(n^2)^d} |\mu|(A_{n^2})
\]
is almost surely finite, and so 
\[
\lim_{n\to\infty}\frac{1}{(n^2)^d} |\mu|(A_{n^2})=0, 
\]
almost surely.
Finally, given $R>0$, let $n\in \N$ be such that
\[
n^2 \leq R \leq (n+1)^2 \,.
\]
Since
\[
A_R \subseteq \{p\in X: n^2<\lVert p\rVert \leq (n+1)^2\} \,\cup \,A_{n^2} \,,
\]
we conclude
\begin{align*}
\frac{1}{R^d} |\mu|(A_R) &\lesssim \frac{1}{R^d} ((n+1)^2)^{d-1}(2n+1) + \frac{1}{R^d}|\mu|(A_{n^2})  \\
&\lesssim \frac{((n+1)^2)^{d-1}}{(n^2)^{d-1}}\frac{2n+1}{n^2} + \frac{|\mu|(A_{n^2)}}{(n^2)^d} \xrightarrow[n\to \infty]{} 0 \,.
\end{align*}

\medskip

For the second limit, notice that
\[
\tilde{A}_R:= \{p\in X:\|p\|> R, \|p+\xi_p\|\leq R\} \subseteq  \{p\in X: R<\|p\|\le R + \|\xi_p\|\}\,. 
\]
Hence
\[
|\mu|(\tilde{A}_R)\leq \sum_{p\in X\setminus X_R} |\mu(p)|\mathbbm{1}_{\{\lVert p\rVert \leq R + \lVert \xi_p\rVert\}} \,.
\]
So
\begin{align*}
\e(|\mu|(\tilde{A}_R)) \leq \sum_{p\in X\setminus X_R} |\mu(p)| \e(\mathbbm{1}_{\{\lVert p\rVert \leq R + \lVert \xi\rVert\}}) 
 &\leq \e\left(\sum_{p\in X_{R+\|\xi\|+1}  \setminus X_R} |\mu(p)|\right) \\[2mm]
 &\lesssim \e((R+\|\xi\|+1)^{d-1}(\|\xi\|+1)) \,.
 \end{align*}
We get
\[
\e\left(\frac{1}{R^d} |\mu|(\tilde{A}_R)\right) \lesssim \frac{1}{R}\bigl(\e(\lVert \xi\Vert^d)+1\bigr) \,.
\]
As before, almost surely
\[
\lim_{n\to\infty} \frac{1}{(n^2)^d} |\mu|(\tilde{A}_{n^2}) = 0 \,.
\]
Finally, given $R>0$, let $n\in \N$ be such that
\[
n^2 \leq R \leq (n+1)^2 \,.
\]
Since
\[
\tilde{A}_R \subseteq \{p\in X: n^2<\lVert p\rVert \leq (n+1)^2\} \,\cup \,\tilde{A}_{(n+1)^2} \,,
\]
we conclude
\begin{align*}
\frac{1}{R^d} |\mu|(\tilde{A}_R) &\lesssim \frac{1}{R^d} ((n+1)^2)^{d-1}(2n+1) + \frac{1}{R^d}|\mu|(\tilde{A}_{(n+1)^2})  \\
&\lesssim \frac{((n+1)^2)^{d-1}}{(n^2)^{d-1}}\frac{2n+1}{n^2} + \frac{((n+1)^2)^d}{(n^2)^d}\frac{|\mu|(\tilde{A}_{(n+1)^2)}}{((n+1)^2)^d} \xrightarrow[n\to \infty]{} 0 \,.
\end{align*}

\end{proof}

\begin{rem}\label{W cs finita}
Note that, in particular, this result implies that, almost surely, \( |\widetilde\mu|(B_R) \) is finite for every $R>0$. Hence, the sums in \eqref{mainresult} are well defined almost surely. Even more: the result implies the almost surely polynomial growth \( |\widetilde\mu|(B_R) \lesssim R^d \). 
\end{rem}

\bdem[Proof of Proposition \ref{key para dividir}]
Note that both sums include the terms corresponding to
\[
\{p\in X: \|p\|\le R, \|p+\xi_p\|\le R\} \,.
\]
Hence
     \begin{align*}
         \left|\frac{1}{R^d} \sum_{p\in X: p+\xi_p\in B_R} e^{-2\pi i \langle p+\xi_p,\la \rangle} \mu(p)- \frac{1}{R^d}\sum_{p\in X_R}  e^{-2\pi i \langle p+\xi_{p}, \la \rangle}\mu(p)\right|\leq (\mbox{I})+(\mbox{II}) \\
     \end{align*}    
     where
       \begin{align*}
         (\mbox{I})&=\left|\frac{1}{R^d} \sum_{\substack{p\in X_R \\  p+\xi_{p} \notin B_R}} e^{-2\pi i \langle p+\xi_p,\la \rangle} \mu(p)\right| 
            \leq\frac{1}{R^d} |\mu| \big(\{p\in X:  \lVert p\rVert \leq R, \lVert p+\xi_{p}\rVert > R\}\big)\\
         \intertext{and}
         (\mbox{II})&=\left|\frac{1}{R^d} \sum_{\substack{p\notin X_R \\ p+\xi_{p} \in B_R}} e^{-2\pi i \langle p+\xi_p,\la \rangle} \mu(p)\right| 
         \leq \frac{1}{R^d} |\mu| \big( \{p\in X:  \lVert p\rVert > R, \lVert p+\xi_{p}\rVert \leq R\} \big)\,.
     \end{align*}
By applying Lemma~\ref{adentroafuera} to each part, we get the desired result.
\edem

\section{The Random Part}\label{random part}

In the previous section, we have reduced the proof of Theorem  \ref{principal1} to prove that almost surely, for every \(\lambda\in\mathbb{R}^d\),
\begin{equation}\label{random remainder}
\lim_{R\to\infty}\frac{1}{m_d(B_R)}\sum_{p\in X_R}\Bigl[e^{-2\pi i\langle p,\lambda\rangle}\mu(p)
\Bigl(e^{-2\pi i\langle \xi_p,\lambda\rangle}-\mathbb{E}\bigl(e^{-2\pi i\langle \xi_p,\lambda\rangle}\bigr)\Bigr)\Bigr]=0 \,.
\end{equation}

For each fixed \(\lambda\), the presence of weights and phase factors leads to a weighted version of the strong law of large numbers. The desired convergence follows directly from Theorem \ref{leyadoc}. The main additional difficulty is to obtain a single event of probability one on which convergence holds for every \(\lambda\in\mathbb R^d\).

\medskip

Recall that $\mu$ is a translation-bounded crystalline measure, with $\mu$ and $\widehat{\mu}$ as in \textnormal{\textbf{%
\textcolor{black}{(}%
\hyperref[eq:N]{N}%
\textcolor{black}{)}}}. Also recall that $ \{\xi_{p}\}_{p\in X}$ is a family of random vectors taking values in $\R^d$, and satisfying condition {\hypersetup{linkcolor=azul}%
\textnormal{\textbf{(\ref{eq:A1})}}} and the tail's control \textbf{(\ref{eq:A2})} by a random vector $\xi$ with finite variance. This section is devoted to proving $\eqref{random remainder}$. Actually, we will prove the following more general result: 

\begin{teo}\label{chuecageneralizada}
Let $\{Y_p(\la)\}_{p\in X,\la\in \R^d}$ be a mean zero family of complex-valued random variables such that
\begin{itemize}
    \item[i.)] $|Y_p(\la)|$ is uniformly bounded over $p\in X$ and $\la\in\R^d$.
    \item[ii.)] For each $\la\in \R^d$, the family $\{Y_p(\la)\}_{p\in X}$ is independent.
    \item[iii.)] For each $p\in X$, $Y_p(\cdot)\in C^1(\R^d)$.
    \item[iv.)] There exists an independent family $\{Z_p\}_{p\in X}$ of real-valued random variables with $\e(Z_p^2)$ uniformly bounded over $p\in X$, such that $\sup_{\la\in \R^d} \lVert \nabla Y_p(\la)\lVert_1 \leq Z_p$, where $\|\cdot\|_1$ denotes the usual $L_1$-norm.
\end{itemize} 
Then, almost surely, for all $\la,\tilde{\la}\in\R^d$
\begin{equation*}
    \lim_{R\to\infty} \frac{1}{R^d} \sum_{p\in X_R} \left[e^{-2\pi i \langle p,\tilde{\la} \rangle} \mu(p)Y_p(\la)\right] = 0 \,.
\end{equation*}

\end{teo}

\medskip

This theorem may be regarded as a \textit{“twisted uniform Law of Large Numbers”}. The uniformity in $\lambda$ and $\tilde{\lambda}$ will be addressed through a union-bound–type argument. From this perspective, the result is closely related to the Wiener–Wintner ergodic theorems \cite{MR1995517}. Various techniques can be employed to tackle problems of this nature; one that has proven to be particularly effective for this specific limit is based on the notion of \textit{affinity}.

\medskip

Before proceeding to the proof of Theorem \ref{chuecageneralizada}, let us show how it can be used to prove \eqref{random remainder}, and in this way to finish the proof of Theorem \ref{principal1}.

\bdem[Proof of \eqref{random remainder} using Theorem \ref{chuecageneralizada}] Take
$$
Y_p(\la):=e^{-2\pi i \langle \xi_p,\la \rangle} - \e(e^{-2\pi i \langle \xi_p,\la \rangle}).
$$
Clearly the family $\{Y_p(\la)\}_{p\in X,\la\in \R^d}$ satisfies the conditions (i)-(iii). On the other hand, with similar calculations to those made in \cite{MR4526342, MR4406130}, we get the estimates

\begin{align*}
    \| \nabla(Y_p) \|_1 \leq   \lVert \nabla \e\bigl(e^{2\pi i \langle \xi_p,\la\rangle}\bigr) \rVert_1+ \lVert \nabla e^{2\pi i \langle \xi_p,\la\rangle} \rVert_1 \leq  2\pi\Big( \e\left(\lVert \xi \rVert_1\right)+  \lVert \xi_p\rVert_1\Big).
    \end{align*}
This proves that (iv) also holds, and therefore \eqref{random remainder} can be deduced as a consequence of Theorem \ref{chuecageneralizada}.
\edem

\begin{rem}
It is worth noting that, although in Theorem \ref{principal1} we impose a stronger condition on the moments of these random variables $\xi_p$, to prove  \eqref{random remainder}, we only need that the second moment be finite. The ``rest'' of the moments were used in Lemma \ref{adentroafuera}.
\end{rem}

\subsection{Proof of Theorem \ref{chuecageneralizada}}

\subsubsection{Affinity and the Hellinger measure}

Let $\nu_1$ and $\nu_2$ be two locally finite positive Radon measures, and let $\sigma$ be another locally finite positive Radon measure such that  $\nu_1$ and $\nu_2$ are absolutely continuous with respect to $\sigma$. As usual, this will be written as  $\nu_1 \ll \sigma$ and $\nu_2 \ll \sigma$. Let
$$
\frac{d\nu_1}{d\sigma} \qquad\mbox{and} \qquad \frac{d\nu_2}{d\sigma} 
$$
denote the Radon–Nikodym derivatives of $\nu_1$ and $\nu_2$ with respect to $\sigma$. Then, the \textit{localized Hellinger measure} is defined as
$$
\rho(\nu_1,\nu_2):=\left[\left(\frac{d\nu_1}{d\sigma}\right)^{1/2}\cdot\quad \left(\frac{d\nu_2}{d\sigma}\right)^{1/2} \right]\ d\sigma.
$$

The total mass of $\rho(\nu_1,\nu_2)$ is often called \textit{affinity} of the measures $\nu_1$ and $\nu_2$. The following three properties of the Hellinger measure can be deduced from its definition without difficulty: 
\begin{itemize}
\item The measure $\rho(\nu_1,\nu_2)$ does not depend on the chosen measure $\sigma$. Typically, the reference measure $\sigma$ is taken as the sum $\nu_1+\nu_2$.
\item It is monotone with respect to each measure.
\item The measures $\nu_1$ and $\nu_2$ are mutually singular if and only if $\rho(\nu_1,\nu_2)=0$. 
\end{itemize}

Another property of the Hellinger measures that we will need to prove Theorem \ref{chuecageneralizada} is the following.

\begin{lem}\label{limsup de afinidades} Let $\nu_t, \gamma_t$ be two families of positive measures over $\mathbb{R}^d$, where $t>0$, such that $\nu_t \rightarrow \nu$ and $\gamma_t \rightarrow \gamma$, in the distributional sense, as $t \rightarrow \infty$. Then, for every non-negative function $f\in C^\infty_c(\R^d)$
$$
 \limsup_{t\rightarrow \infty}\int_{\R^d} f(x)\, d\rho\left(\nu_t, \gamma_t\right) \leq \int_{\R^d} f(x)\, d\rho\left(\nu, \gamma\right).
$$
\end{lem}

The interested reader is referred to \cite[Thm. 5]{MR472749} or \cite[Appendix A2]{MR4526342} for a proof of this result.

\bigskip

\subsubsection{Connection between the Hellinger measure and our problem}\label{connection}

The connection between Theorem \ref{chuecageneralizada} and the notion of affinity is another critical point in our argument, as this is precisely where earlier works make the strongest use of the separation of the set $X$. In our case, although $X$ is locally finite and the measure $\mu$ is translation-bounded, the assignment of weights $\mu(p)$ to each delta $\delta_p$ allows $X$ to fail to be uniformly separated; indeed, it may not even be a finite union of uniformly separated sets, as happens in the examples of crystalline measures given by Kolountzakis \cite{MR3448912} and Meyer \cite{MR3482845}. This section shows how to overcome this difficulty and, using the notion of Hellinger measure, essentially proves Theorem \ref{chuecageneralizada}. We say essentially, because the proof of Proposition \ref{difracaleatoria} below, which is rather long and more technical, is left to the next subsection.  

\medskip

To begin with, let us define

\begin{equation}\label{eseerre}
S_R(\tilde{\lambda},\lambda):=\frac{1}{m_d(B_R)}\sum_{p\in X_R}e^{2\pi i\langle p,\tilde{\lambda}\rangle}\mu(p)\,Y_p(\la) \,.
\end{equation}

Fix \(\lambda,\tilde{\lambda}\in\mathbb{R}^d\), and  choose \(\gamma_0\in X^*\). Since $X^*$ is locally finite, there is some $\eps>0$ such that 
\[
\lVert\gamma_0-\gamma\rVert\geq \eps, \quad \text{ for all } \gamma\in X^*\setminus\{\gamma_0\} \,.
\]
Let \(f:=f_{\tilde{\lambda},\gamma_0}\in C_c^\infty(\R^d)\) such that
\begin{equation*}
    f(\w-\tilde{\la}) = \begin{cases} \ds \frac{1}{A(\gamma_0)} & \mbox{if $\ \w=0$}\\[3mm]
    0 & \mbox{if $\|\w\|\geq \eps$}
    \end{cases} .
\end{equation*}
In particular, $f(\gamma-\gamma_0-\tilde{\la}) = 0$ for every $\gamma\in X^*\setminus\{\gamma_0\}$. Then by the Poisson summation formula associated to the crystalline measure $\mu$
\begin{align*}
\sum_{q\in X}\mu(q)\, \Big(e^{-2\pi i\langle \tilde{\lambda}+\gamma_0, q\rangle}\widehat{f}(q-p)\Big)
&=
\sum_{\gamma\in X^*}A(\gamma) \Big(f(\gamma-\gamma_0-\tilde{\lambda})e^{2\pi i\langle \gamma-\gamma_0-\tilde{\la},p\rangle}\Big)\\
&=e^{-2\pi i\langle \tilde{\lambda},p\rangle} \,.
\end{align*}
Since $\overline{\widehat{f}(\w)}=\widehat{\overline{f}}(-\w)$, we have 
\begin{align*}
S_R(\tilde{\lambda},\lambda)&=\frac{1}{m_d(B_R)}\sum_{p\in X_R}\left[\left(\sum_{q\in X}\overline{\mu(q)e^{-2\pi i\langle \tilde{\lambda}+\gamma_0,q\rangle}}\,\widehat{\overline{f}}(p-q)\right)
\mu(p)\, Y_p(\la)\right]\\
&=\frac{1}{m_d(B_R)}\sum_{p\in X_R}\sum_{q\in X}\Big[\overline{\mu(q)}e^{2\pi i\langle \tilde{\lambda}+\gamma_0,q\rangle}\Big]\cdot
\Big[\mu(p)\, Y_p(\la) \Big]\,\widehat{\overline{f}}(p-q) \,.
\end{align*}

So, if we define the translation-bounded measures
\begin{eqnarray*}
 \mu_{\alpha}=\sum_{p\in X} e^{-2\pi i \langle \alpha, p\rangle}\mu(p)\,\delta_p  &\mbox{and}& \nu_{\la}=\sum_{p\in X} Y_p(\la) \mu(p)\, \delta_p, 
\end{eqnarray*}
and we define $\mu_\alpha(p):=e^{-2\pi i \langle \alpha, p\rangle}\mu(p)$ and $\nu_\la(p):= Y_p(\la) \mu(p)$, then 
\begin{align*}
S_R(\tilde{\lambda},\lambda)&=\frac{1}{m_d(B_R)}\sum_{p\in X_R}\sum_{q\in X}\Big( \overline{\mu_{\tilde{\lambda}+\gamma_0}(q)}\cdot \nu_{\la}(p)\Big) \,\widehat{\overline{f}}(p-q)\\
&=\frac{1}{m_d(B_R)}\left(\sum_{p,q\in X_R}\Big( \overline{\mu_{\tilde{\lambda}+\gamma_0}(q)}\cdot \nu_{\la}(p)\Big)\right. 
+\left.\sum_{p\in X_R}\sum_{q\in X\setminus X_R}\Big( \overline{\mu_{\tilde{\lambda}+\gamma_0}(q)}\cdot \nu_{\la}(p)\Big) \right) \,\widehat{\overline{f}}(p-q).
\end{align*}

By Proposition \ref{sevuelveausar}, we already know that
$$
\lim_{R\to\infty} \frac{1}{m_d(B_R)} \sum_{p\in X_R}\sum_{q\in X\setminus X_R}\Big( \overline{\mu_{\tilde{\lambda}+\gamma_0}(q)}\cdot \nu_{\la}(p)\Big) \,\widehat{\overline{f}}(p-q)=0.
$$
So, to prove Theorem \ref{chuecageneralizada}, we need to show that (almost surely)
\begin{equation}\label{eq ya casi}
\lim_{R\to\infty} \frac{1}{m_d(B_R)} \sum_{p\in X_R}\sum_{q\in X_R}\Big( \overline{\mu_{\tilde{\lambda}+\gamma_0}(q)}\cdot \nu_{\la}(p)\Big)\,\widehat{\overline{f}}(p-q) = 0 \,.
\end{equation}

Recall that for a purely atomic measure $\mu_0$ whose support is contained in $X$
\begin{equation*}
    \gamma_{\mu_0}^{(R)} = \frac{1}{m_d(B_R)} \sum_{p,q\in X_R} \mu_0(p)\overline{\mu_0(q)} \delta_{p-q} \,,
\end{equation*}
then
\begin{align*}
\widehat{\gamma_{\mu_{\tilde{\lambda}+\gamma_0}}^{(R)}}&= \frac{1}{m_d(B_R)}  \left| \sum_{p\in X_R}  \mu_{\tilde{\lambda}+\gamma_0}(p) e^{-2\pi i\langle x, p\rangle} \right|^2 \, dx\\
\widehat{\gamma_{\nu_{\la}}^{(R)}}&= \frac{1}{m_d(B_R)}  \left| \sum_{p\in X_R}  \nu_{\lambda}(p) e^{-2\pi i\langle x, p\rangle}  \right|^2 \, dx\,.
\end{align*}

Hence, an easy computation shows that
\begin{equation}\label{eq ya casi casi}
 \frac{1}{m_d(B_R)}\left|\sum_{p\in X_R}\sum_{q\in X_R}\Big( \overline{\mu_{\tilde{\lambda}+\gamma_0}(q)}\cdot \nu_{\la}(p)\Big)\,\widehat{\overline{f}}(p-q)\right|\leq  \int_{\R^d} |\overline{f}|(x)\, d\rho\bigl(\widehat{\gamma_{ \mu_{\tilde{\lambda}+\gamma_0}}^{(R)}},\widehat{\gamma_{ \nu_{\la}}^{(R)}}\bigr).
\end{equation}

By Corollary \ref{diffracdeterm}, the measures  $\widehat{\gamma_{\mu_{\tilde{\lambda}+\gamma_0}}^{(R)}}$ converge, in the sense of distributions, to the purely atomic measure 
\[
\widehat{\gamma_{\mu_{\tilde{\lambda}+\gamma_0}}}=\sum_{u\in X^*}|A(u)|^2\,\delta_{(u-(\tilde{\lambda}+\gamma_0))} \,.
\]

On the other hand, we have the following result. 

\begin{pro}\label{difracaleatoria}
The measure $\gamma_{\nu_{\la}}^{(R)}$ 
can be split as
$$
\gamma_{\nu_{\la}}^{(R)}=\gamma_{\la,R}+ C_{\la,R} \delta_0,
$$
where $C_{\la,R}$  is bounded independently of $\la$ and $R$, and $\gamma_{\la,R}$ is a random measure such that almost surely, for every $\la\in \R^d$:
\begin{equation}\label{almostsurelyvag}
    \lim_{R\to\infty}\gamma_{\la,R}= 0 \,,
\end{equation}
where the convergence is in the weak topology. 
\end{pro}

\medskip

Using this result, we can get $\eqref{eq ya casi}$ and in consequence, get Theorem \ref{chuecageneralizada}. Let us work in the full probability set where the limit $\eqref{almostsurelyvag}$ holds. Now, it is enough to prove that for any increasing sequence of positive real numbers $\{R_n\}_{n\in\N}$, there exist a subsequence $\{R_{n_k}\}_{k\in\N}$ such that 
\begin{equation*}
\lim_{k\to\infty} \frac{1}{m_d(B_{R_{n_k}})} \sum_{p\in X_{R_{n_k}}}\sum_{q\in X_{R_{n_k}}}\Big( \overline{\mu_{\tilde{\lambda}+\gamma_0}(q)}\cdot \nu_{\la}(p)\Big)\,\widehat{\overline{f}}(p-q) = 0 \,.
\end{equation*}
So, given such a sequence $\{R_n\}_{n\in\N}$, by compactness, there exists a subsequence $S=\{R_{n_k}\}_{k\in\N}$ such that $C_{\la,R_{n_k}}\to C_{\la,S}$ as $k\to\infty$. Then, by Proposition \ref{difracaleatoria}
$$
\gamma_{\nu_{\la}}^{(R_{n_k})}\xrightarrow[k\to\infty]{\mbox{\small weakly}}C_{\la,S} \delta_0.
$$ 
Since the translation-bounded norms of the measures $\gamma_{\nu_{\la}}^{(R)}$ are uniformly bounded, the convergence is also in the distributional sense (Lemma \ref{vagimpliesdist}).
Therefore, the measures 
$$
\widehat{{\gamma}_{\nu_{\la}}^{(R_{n_k})}}\xrightarrow[k\to\infty]{\mbox{\small in the distributional sense}}C_{\la,S} dx.
$$ 
Therefore, we can conclude the proof of Theorem \ref{chuecageneralizada} as follows. By equation \eqref{eq ya casi casi}, we obtain that 
\begin{align*}
\limsup_{k\to\infty} \frac{1}{m_d(B_{R_{n_k}})}\left|\sum_{p,q\in X_{R_{n_k}}}\Big( \overline{\mu_{\tilde{\lambda}+\gamma_0}(q)}\cdot \nu_{\la}(p)\Big)\,\widehat{\overline{f}}(p-q)\right| \leq\limsup_{k\to\infty} \int_{\R^d} |\overline{f}|(x)\, d\rho\bigl(\widehat{\gamma_{ \mu_{\tilde{\lambda}+\gamma_0}}^{(R_{n_k})}},\widehat{\gamma_{ \nu_{\la}}^{(R_{n_k})}}\bigr).
\end{align*}

Recall that, given a function $f\in\schw$, there exists $g\in\schw$ such that $|f|\leq g$. In our case, we need a majorant for $\overline{f}$, which has compact support. Hence, the function $g$ can be chosen with compact support. So, using this and Lemma \ref{limsup de afinidades}, we get that 
\begin{align*}
\limsup_{k\to\infty} \int_{\R^d} |\overline{f}|(x)\, d\rho\Bigl(\widehat{\gamma_{ \mu_{\tilde{\lambda}+\gamma_0}}^{(R_{n_k})}},\widehat{\gamma_{ \nu_{\la}}^{(R_{n_k})}}\Big) &\leq
\limsup_{k\to\infty} \int_{\R^d} g(x)\, d\rho\Bigl(\widehat{\gamma_{ \mu_{\tilde{\lambda}+\gamma_0}}^{(R_{n_k})}},\widehat{\gamma_{ \nu_{\la}}^{(R_{n_k})}}\Big)\\& 
\leq\int_{\R^d} g(x)\, d\rho\Bigl(\widehat{\gamma_{\mu_{\tilde{\lambda}+\gamma_0}}},C_{\la,S} dx\Big) = 0, 
\end{align*}
where for the last identity we have used that the measures $C_{\la,S} dx$ and $\widehat{\gamma_{\mu_{\tilde{\lambda}+\gamma_0}}}$ are mutually singular. It only remains to prove Proposition \ref{difracaleatoria}. The next subsection is devoted to this aim.

\subsubsection{Proof of Proposition \ref{difracaleatoria}}\label{sequential}

Note that, by the independence of the random variables $\nu_\la(p)$ and  $\nu_\la(q)$ for each $p\neq q$, and the zero means $\e(\nu_\la(p))$, the measure $\gamma_{\nu_\la}^{(R)}$ can be rewritten as follows:
\begin{align*}
    \gamma_{\nu_\la}^{(R)}=\frac{1}{m_d(B_{R})}\left( \sum_{p,q\in X_{R}} \left[\nu_\la(p) \overline{\nu_\la(q)}-\e\left(\nu_\la(p) \overline{\nu_\la(q)}\right)\right]\delta_{(p-q)} + \delta_0\sum_{p\in X_{R}} \e\left(|\nu_\la(p)|^2\right)\right)
\end{align*}
and 
$$
\left(\frac{1}{m_d(B_{R})} \sum_{p\in X_{R}} \e\left(|\nu_\la(p)|^2\right)\right) \delta_0=\left(\frac{1}{m_d(B_{R})} \sum_{p\in X_R}|\mu(p)|^2 \mbox{Var}(Y_p(\la)) \right) \delta_0.
$$
Since the quantity 
\begin{equation}\label{quantity}
\frac{1}{m_d(B_{R})} \sum_{p\in X_R}|\mu(p)|^2 \mbox{Var}(Y_p(\la)) 
\end{equation}
can be bounded independently of $\la$ and $R$, the proof of Proposition \ref{difracaleatoria} would be complete if we can prove the following result.
\begin{prop}\label{delKalasf}
Almost surely, for every \(\lambda\in\mathbb{R}^d\) 
\begin{equation}\label{eqdelKalasf}
\lim_{R\to\infty}\frac{1}{m_d(B_R)}\sum_{p,q\in X_R}
\Bigl[\nu_\lambda(p)\overline{\nu_\lambda(q)}-\mathbb{E}\bigl(\nu_\lambda(p)\overline{\nu_\lambda(q)}\bigr)\Bigr]f(p-q)=0 \,,
\end{equation}
for every $f\in C_c(\R^d)$.
\end{prop}

The proof follows the general scheme of \cite{MR4526342}, and it is based on the next lemmas. However, to include the weights in each atom, we need to use different tools, like some weighted versions of the Law of Large Numbers. 

\begin{lem}\label{Kacotado1}
Let \(K\subseteq \mathbb{R}^d\) be a bounded set. Given $\la\in \R^d$, the following limit holds almost surely: 
\[
\lim_{R\to\infty}\frac{1}{m_d(B_R)}\sum_{p\in X_R}\left(\sum_{q\in X_R\cap(K+p)}\nu_\lambda(p)\overline{\nu_\lambda(q)}
-\mathbb{E}\bigl(\nu_\lambda(p)\overline{\nu_\lambda(q)}\bigr)\right)=0 \,.
\]
\end{lem}

\begin{proof}
Recall that for each \(\lambda\in\mathbb{R}^d\), the random variables $\nu_\lambda(p)=Y_p(\la) \mu(p)$, with $p\in X$, constitute an independent family of random variables with zero mean. So we can separate
    \begin{align*}
        &\frac{1}{m_d(B_R)} \sum_{p\in X_R}\left(\sum_{q\in X_R\cap (K+p)} \nu_\la(p) \overline{\nu_\la(q)}-\e\bigl(\nu_\la(p) \overline{\nu_\la(q)}\bigr)\right)  = \\
        &\frac{1}{m_d(B_R)} \sum_{p\in X_R}\left(\sum_{q\in X_R\cap ((K+p)\setminus\{p\})} \nu_\la(p) \overline{\nu_\la(q)} \right) + \mathbbm{1}_K(0)\frac{1}{m_d(B_{R})}\sum_{p\in X_{R}} \left(|\nu_\la(p)|^2 - \e\left(|\nu_\la(p)|^2\right)\right) \,.
\end{align*}

\noindent For the second term, let 
    \begin{equation*}
        \mu^2:= \sum_{p\in X} \mu(p)^2\delta_p \,.
\end{equation*}
If 
\begin{equation*}
    |\mu^2|(X)<\infty 
\end{equation*}
then
\begin{align*}
    &\left|\frac{1}{m_d(B_{R})}\sum_{p\in X_{R}} \left(|\nu_\la(p)|^2 - \e\left(|\nu_\la(p)|^2\right)\right)\right| \lesssim \frac{|\mu^2|(B_R)}{m_d(B_R)} \leq  \frac{|\mu^2|(X)}{m_d(B_R)} \xrightarrow[R\to\infty]{} 0 \,.
\end{align*}
If
\begin{equation*}
    |\mu^2|(X) = \infty \,,
\end{equation*}
let $M>0$ be such that $|\mu(p)|\leq M$, for each $p\in X$. Notice that
    \begin{equation}\label{comparacion de medidas}
        \frac{|\mu^2|(B_R)}{m_d(B_R)} \leq \frac{M|\mu|(B_R)}{m_d(B_R)} \leq M' \,,
    \end{equation}
then
     \begin{align*}
          &\left|\frac{1}{m_d(B_{R})}\sum_{p\in X_{R}} \left(|\nu_\la(p)|^2 - \e\left(|\nu_\la(p)|^2\right)\right)\right| \leq M'\left|\frac{1}{|\mu^2|(B_R)}\sum_{p\in X_{R}} |\mu(p)|^2\left(|Y_p(\la)|^2 - \e\left(|Y_p(\la)|^2\right)\right)\right|\,.
     \end{align*}
So we can apply Theorem \ref{leyadoc} to prove that this term tends to zero as $R\to\infty$. 
For the first term, remember that $|\mu|(X) = \infty$ (Corollary \ref{muinfinita}). Let
     \begin{equation*}
         W_p := Y_p(\la) \sum_{q\in X\cap ((K+p)\setminus\{p\})} \overline{\mu(q)  Y_q(\la)} \,.
     \end{equation*}
Then, using \eqref{comparacion de medidas}, we can make the estimate 
     \begin{equation*}
         \left|\frac{1}{m_d(B_R)} \sum_{p\in X_R}\left(\sum_{q\in X\cap ((K+p)\setminus\{p\})} \nu_\la(p) \overline{\nu_\la(q)} \right)\right| \lesssim \left|\frac{1}{|\mu|(B_R)} \sum_{p\in X_R} \mu(p)W_p\right| \,
     \end{equation*}
and then, by Theorem \ref{leyadoc} we get that
    \begin{equation*}
        \lim_{R\to\infty} \frac{1}{m_d(B_R)} \sum_{p\in X_R}\left(\sum_{q\in X\cap ((K+p)\setminus\{p\})} \nu_\la(p) \overline{\nu_\la(q)} \right) = 0 \,,
    \end{equation*}
almost surely. Now, if $R_0\geq 1$ is such that $K\subseteq B_{R_0}$, and $\lVert p\rVert \leq R$, then
    \begin{equation*}
        X \cap (K+p) = X_{R+R_0} \cap (K+p) \,,
    \end{equation*}
hence
    \begin{align*}
        &\frac{1}{m_d(B_R)} \sum_{p\in X_R}\left(\sum_{q\in X_R\cap ((K+p)\setminus\{p\})} \nu_\la(p) \overline{\nu_\la(q)} \right) = \\
        &\frac{1}{m_d(B_R)} \sum_{p\in X_R}\left(\sum_{q\in X\cap ((K+p)\setminus\{p\})} \nu_\la(p) \overline{\nu_\la(q)} \right) -\\
        &\frac{1}{m_d(B_R)} \sum_{p\in X_R}\left(\sum_{q\in (X_{R+R_0}\setminus X_R)\cap ((K+p)\setminus\{p\})} \nu_\la(p) \overline{\nu_\la(q)} \right) \,.
    \end{align*}
To prove that the last sum goes to $0$, notice that $K+p\subseteq B_R$ provided that $\lVert p \rVert\leq R-R_0$. Then
    \begin{align*}
        &\left|\frac{1}{m_d(B_R)} \sum_{p\in X_R}\left(\sum_{q\in (X_{R+R_0}\setminus X_R)\cap ((K+p)\setminus\{p\})} \nu_\la(p) \overline{\nu_\la(q)} \right)\right| =\\[2mm]
        &\left|\frac{1}{m_d(B_R)} \sum_{p\in X_R\setminus X_{R-R_0}}\left(\sum_{q\in (X_{R+R_0}\setminus X_R)\cap ((K+p)\setminus\{p\})} \nu_\la(p) \overline{\nu_\la(q)} \right)\right| \lesssim \\[2mm]
        &\frac{1}{m_d(B_R)} \sum_{p\in X_R\setminus X_{R-R_0}}|\mu(p)|\left(\sum_{q\in (X_{R+R_0}\setminus X_R)\cap ((K+p)\setminus\{p\})}  |\mu(q)| \right) \lesssim \\[2mm]
        &\frac{R_0^{d+1} R^{d-1}}{m_d(B_R)}\xrightarrow[\quad R\to \infty\quad ]{} 0\,
    \end{align*}
and this finishes the proof.
\edem

Now, we will prove that we can take a set of probability one that is common for every $\la\in\R^d$. The proof follows a union-bound-type argument similar to the one used in  \cite{MR4406130} and \cite{MR4526342}.

\begin{lem}\label{Kacotad2}
Let \(K\subseteq \mathbb{R}^d\) be a bounded set. Then, almost surely for every \(\lambda\in\mathbb{R}^d\),
\begin{equation}\label{paratodolambda}
\lim_{R\to\infty}\frac{1}{m_d(B_R)}\sum_{p\in X_R}\left(\sum_{q\in X_R\cap(K+p)}\nu_\lambda(p)\overline{\nu_\lambda(q)}
-\mathbb{E}\bigl(\nu_\lambda(p)\overline{\nu_\lambda(q)}\bigr)\right)=0 \,.
\end{equation}
\end{lem}

Before proceeding to the proof, we state, without demonstration, the following simple but useful result.

\begin{lem}\label{aux1}
    Let $(\Omega,\efe,\pe)$ be a probability space and $g:\Omega \times \R^d \to \ce$ a function such that

 \begin{itemize}
 
 \item For each $\omega\in \Omega$, $g(\omega,\cdot)\in C^1(\R^d)$. 

 \item There exists $\theta:\Omega\to\R^d$ with $\e(\theta)<\infty$, and for every $\omega\in \Omega$ and every $1\leq j\leq d$: $$\sup_{\la\in\R^d} \left| \frac{\partial g(\omega,\la)}{\partial \la_j}\right| \leq \theta(\omega) \,.$$
 \end{itemize}

Then for each $1\leq j \leq d$: 
\begin{equation*}
    \frac{\partial \e(g(\cdot,\la))}{\partial \la_j} = \e\left(\frac{\partial g(\cdot,\la)}{\partial \la_j}\right) \in C(\R^d) \,.
\end{equation*}
\end{lem}

\medskip

\bdem[Proof of Lemma \ref{Kacotad2}] Let
\begin{align*}
    F_{R}(\la)&:= \frac{1}{R^d} \sum_{p\in X_R}\sum_{q\in X_R\cap ((K+p)\setminus\{p\})} \nu_\la(p)\overline{\nu_\la(q)} \,,\\
    G_R(\la)&:= \frac{1}{R^d}\sum_{p\in X_{R}} \left(|\nu_\la(p)|^2 - \e\left(|\nu_\la(p)|^2\right)\right) \,.
\end{align*}
We claim that, almost surely 
\begin{align*}
    \sup_{R\geq 1}\sup_{\la\in \R^d} \lVert \nabla F_R(\la)\rVert &< \infty  \,, \\ 
    \sup_{R\geq 1}\sup_{\la\in \R^d} \lVert \nabla G_R(\la)\rVert &< \infty \,.
\end{align*}

\medskip

On the one hand,
\begin{align*}
    \lVert \nabla F_R(\lambda) \rVert &\leq \frac{1}{R^d} \sum_{p \in X_R} \left( |\mu(p)| \lVert \nabla Y_p(\lambda) \rVert \sum_{q \in X_R \cap ((K+p) \setminus \{p\})} |\mu(q)| |Y_q(\lambda)| \right) \\
    &\quad + \frac{1}{R^d} \sum_{p \in X_R} \left( |\mu(p)| |Y_p(\lambda)| \sum_{q \in X_R \cap ((K+p) \setminus \{p\})} |\mu(q)| \lVert \nabla Y_q(\lambda) \rVert \right) \\
    &\lesssim \frac{1}{R^d} \sum_{p \in X_R} \left( |\mu(p)| Z_p \sum_{q \in X_R \cap ((K+p) \setminus \{p\})} |\mu(q)| \right) \\
    &\quad + \frac{1}{R^d} \sum_{p \in X_R} \left( |\mu(p)| \sum_{q \in X_R \cap ((K+p) \setminus \{p\})} |\mu(q)| Z_q \right).
    \end{align*}

For the first sum: use the translation boundedness and Theorem \ref{leyadoc}. For the second sum 
\begin{align*}
    &\frac{1}{R^d}\sum_{p\in X_R}\left(|\mu(p)|\sum_{q\in X_R\cap ((K+p)\setminus\{p\})} |\mu(q)|\,Z_q \right)\leq \\
    &\frac{1}{R^d}\sum_{p\in X_R}\left(|\mu(p)|\sum_{q\in X\cap (K+p)\setminus\{p\}} |\mu(q)|\, (Z_q - \e(Z_q)) \right) + \\
    &\left(1+\sup_{p'\in X}\e(Z_{p'}^2)\right)\frac{1}{R^d}\sum_{p\in X_R}\left(|\mu(p)|\sum_{q\in X\cap (K+p)\setminus\{p\}} |\mu(q)| \right).
\end{align*}
Therefore, using that $\mu$ is translation-bounded and Theorem \ref{leyadoc}, we obtain that this sum is also bounded. In consequence
\begin{equation*}
    \sup_{R\geq 1}\sup_{\la\in \R^d} \lVert \nabla F_R(\la)\rVert < \infty
\end{equation*}
almost surely. On the other hand, in the case of $G_R(\la)$, note that 
\begin{equation}\label{grad}
      \nabla\e\left(|\nu_\la(p)|^2\right) = \e\left(\nabla|\nu_\la(p)|^2\right) \,,
\end{equation}
Indeed, this follows from Lemma \ref{aux1} taking $g(\cdot,\la):=|\nu_\la(p)|^2$. This function satisfies the hypothesis of this lemma because
\begin{equation*}
    \left|\frac{\partial g(\cdot,\la)}{\partial \la_j}\right| = 2|\mu(p)|^2 \,\re\left(\frac{\partial Y_p(\la)}{\partial \la_j}\, \overline{Y_p(\la)}\right) \lesssim  Z_p \,.
\end{equation*}
Therefore, 
%
%
\begin{align*}
    \lVert \nabla G_R(\la)\lVert &\leq \frac{1}{R^d}\sum_{p\in X_R} \left\lVert \nabla\left(|\nu_\la(p)|^2\right)\right\rVert + \frac{1}{R^d}\sum_{p\in X_R} \left\lVert \e\left(\nabla|\nu_\la(p)|^2\right)\right\rVert \\
    &\lesssim \frac{1}{R^d} \sum_{p\in X_R} |\mu(p)|^2  \lVert \nabla Y_p(\la)\rVert + \frac{1}{R^d}\sum_{p\in X_{R}} \e\left\lVert \nabla|\nu_\la(p)|^2\right\lVert_1 \\
    &\lesssim \frac{1}{R^d}\sum_{p\in X_{R}} |\mu(p)|^2 Z_p  + \frac{1}{R^d}\sum_{p\in X_{R}} |\mu(p)|^2\e\left(Z_p \right) \,.
\end{align*}
Hence, we can use the estimate $|\mu(p)|^2\leq M |\mu(p)|$, the translation boundedness of $\mu$, and Theorem \ref{leyadoc} to conclude that $ \lVert \nabla G_R(\la)\lVert$ is also uniformly bounded.

\medskip

By Lemma \ref{Kacotado1} and the fact that $\mathbb{Q}^d$ is countable, we have that $\eqref{paratodolambda}$ almost surely holds for all $\la'\in \mathbb{Q}^d$. Finally, for $\la\in \R^d\setminus \mathbb{Q}^d$, the mean value estimate gives
$$
\limsup_{R\to\infty}\bigl(|F_R(\lambda)|+|G_R(\lambda)|\bigr)
\lesssim \|\lambda-\lambda'\|,
\qquad \lambda'\in\mathbb Q^d.
$$

Letting \(\lambda'\to\lambda\) proves the conclusion.

\end{proof}

Now we can proceed to the proof of Proposition \ref{delKalasf} to obtain Theorem \ref{chuecageneralizada} and thus conclude this section.

\begin{proof}[Proof of Proposition \ref{delKalasf}]

Let \(f\) be a simple function with compact support. By Lemma \ref{Kacotad2}, the limit $\eqref{eqdelKalasf}$ holds almost surely for every $\la\in\R^d$ for $f$. Then by standard arguments we arrive to any continuous function with compact support (using again the translation boundedness). Finally, recall that there is a countable family $\Phi\subseteq C_c(\R^d)$, such that for all $f\in C_c(\R^d)$, there is a compact set $K_f$ and a sequence $\{\phi_n\}_{n\in \N}\subseteq \Phi$ satisfying 
\[
\operatorname{supp}(\phi_n) \,\cup \, \operatorname{supp}(f) \subseteq K_f\,, \text{ for all } n
\]
and
\[
\lim_{n\to\infty} \left(\sup_{x\in K_f} |f(x)-\phi_n(x)|\right) = 0 \,.
\]
With this, we get $\eqref{eqdelKalasf}$ uniformly for all the set $ C_c(\R^d)$.

\end{proof}

\section{General case}\label{generalization}

As stated in the Introduction, Theorem~\ref{principal1} addresses random perturbations of the support points. But, as we have explained, in the setting of this paper it is natural to allow perturbations of the atomic masses. This section shows how the proof of Theorem~\ref{principal1} can be adapted to obtain Theorem~\ref{principal2}. So, for each $p\in X$ we fix $n(p)\in\mathbb{N}$ and define the random perturbation

$$
\widetilde{\mu} :=\sum_{p\in X}\sum_{k=1}^{n(p)} \big(\alpha_{p,k}\,\mu(p)\big) \, \delta_{p + \xi_{p,k}},
$$
and recall the hypotheses \textnormal{\textbf{(\ref{eq:H1})}-\textbf{(\ref{eq:H6})}} on the random vectors $\xi_{p,k}$ and the complex random variables $\alpha_{p,k}$.

\medskip

The idea, as in the proof of Theorem \ref{principal1}, is to divide the problem in  three parts: structured part, almost periodic part, and random part. The first step is to prove the following generalization of Lemma \ref{adentroafuera}. 

\begin{lem}\label{adentroafuerageneralized}
    Almost surely,
    \[
    \lim_{R\to\infty} \frac{1}{R^d} \sum_{p\in X_R} |\mu(p)|\sum_{\substack {k\in \{1,\ldots,n(p)\}\\p+\xi_{p,k}\notin B_R} } | \alpha_{p,k}| = 0
    \]
    and
    \[
    \lim_{R\to\infty} \frac{1}{R^d} \sum_{p\notin X_R} |\mu(p)|\sum_{\substack {k\in \{1,\ldots,n(p)\}\\p+\xi_{p,k}\in B_R} } | \alpha_{p,k}| = 0 \,.
    \]

\end{lem}

\begin{proof}[Sketch of the proof.]
    Notice that 
    \[
    \sum_{p\in X_R} |\mu(p)|\sum_{\substack {k\in \{1,\ldots,n(p)\}\\p+\xi_{p,k}\notin B_R} } | \alpha_{p,k}| = \sum_{p\in X_R} |\mu(p)|\sum_{k=1}^{n(p)} |\alpha_{p,k}|\,\mathbbm{1}_{\{p+\xi_{p,k}\notin B_R\}}
    \]
    and so
    \begin{align*}
    \e\left(\sum_{p\in X_R} |\mu(p)|\sum_{\substack {k\in \{1,\ldots,n(p)\}\\p+\xi_{p,k}\notin B_R} } | \alpha_{p,k}|\right) &= \sum_{p\in X_R} |\mu(p)|\sum_{k=1}^{n(p)} \e(|\alpha_{p,k}|)\,\e(\mathbbm{1}_{p+\xi_{p,k}\notin B_R}) \\
    &\leq \sum_{p\in X_R} |\mu(p)|\sum_{k=1}^{n(p)} \e(|\alpha_{p,k}|)\,\e(\mathbbm{1}_{\|p\|>R-\|\xi_{p,k}\|}) \\[2mm]
    &\leq \sum_{p\in X_R } |\mu(p)|\e(\mathbbm{1}_{\|p\|>R-\|\xi\|}) \e\left(\sum_{k=1}^{n(p)} |\alpha_{p,k}|\right) \\[2mm]
    &\lesssim \sum_{p\in X_R } |\mu(p)|\e(\mathbbm{1}_{\|p\|>R-\|\xi\|})  \,.
    \end{align*}
    Repeating the argument of Lemma \ref{adentroafuera}, we obtain almost surely
    \[
    \lim_{n\to\infty} \frac{1}{(n^2)^d}\sum_{p\in X_{n^2}} |\mu(p)|\sum_{\substack {k\in \{1,\ldots,n(p)\}\\p+\xi_{p,k}\notin B_{n^2}} } | \alpha_{p,k}| = 0 \,.
    \]
    Finally, given $R>0$, let $n\in \N$ be such that
\[
n^2 \leq R \leq (n+1)^2 \,.
\]
Since
\[
\sum_{p\in X_R} |\mu(p)|\sum_{\substack {k\in \{1,\ldots,n(p)\}\\p+\xi_{p,k}\notin B_R} } | \alpha_{p,k}| \lesssim |\mu|\left(X_{(n+1)^2}\setminus X_{n^2}\right) + \sum_{p\in X_{n^2}} |\mu(p)|\sum_{\substack {k\in \{1,\ldots,n(p)\}\\p+\xi_{p,k}\notin B_{n^2}} } | \alpha_{p,k}| \,,
\]
we conclude in the same way as before. The second limit follows analogously.
\end{proof}

\begin{rem}  
As in the case of Lemma \ref{adentroafuera}, this lemma ensures that the sums in \eqref{mainresultgeneralized} are almost surely well-defined.
\end{rem}

\medskip

Recall that we want to prove 
$$
       \lim_{R\to\infty} \frac{1}{m_d(B_R)} \int_{B_R}e^{-2\pi i \langle x,\la \rangle} \,d\widetilde\mu(x) =  \left\{ \begin{array}{lcc}
             w(\lambda)A(\la) &   if  & \lambda\in X^* \\
             \\ 0 &  if  & \lambda\notin X^*
             \end{array},
   \right.
$$
and the integral in the left hand side can be rewritten as
\begin{equation}\label{integral como suma general}
\sum_{p\in X} \ \sum_{\substack {k\in \{1,\ldots,n(p)\}\\ p+\xi_{p,k}\in B_R} }e^{-2\pi i\langle p+\xi_{p,k},\lambda\rangle}\mu(p)\alpha_{p,k}.
\end{equation}
By Lemma \ref{adentroafuerageneralized} we get, as in the proof of Proposition  \ref{key para dividir}, that
  \begin{equation*}
         \frac{1}{R^d}\left|\sum_{p\in X} \sum_{\substack {k\in \{1,\ldots,n(p)\}\\ p+\xi_{p,k}\in B_R} }e^{-2\pi i\langle p+\xi_{p,k},\lambda\rangle}\mu(p)\alpha_{p,k} - \sum_{p\in X_R} \mu(p)\sum_{k=1}^{n(p)} \alpha_{p,k}\, e^{-2\pi i \langle p+\xi_{p,k}, \la \rangle}\right| \xrightarrow[R\to\infty]{} 0 .
    \end{equation*}
uniformly on $\la\in\R^d$. Hence, it is enough to prove that 
$$
       \lim_{R\to\infty} \frac{1}{m_d(B_R)} \sum_{p\in X_R} \mu(p)\left(\sum_{k=1}^{n(p)} \alpha_{p,k}\, e^{-2\pi i \langle p+\xi_{p,k}, \la \rangle}\right)
=  \left\{ \begin{array}{lcc}
             w(\lambda)A(\la) &   if  & \lambda\in X^* \\
             \\ 0 &  if  & \lambda\notin X^*
             \end{array},
   \right.
$$

instead of the sums in \eqref{integral como suma general}. This opens the possibility to split the problem in three parts: a structured part, and almost periodic part, and a random part.  To simplify the expressions of each of these parts, let us introduce the notation:
$$
\Xi_p(\la)=\sum_{k=1}^{n(p)} \alpha_{p,k}\, e^{-2\pi i \langle \xi_{p,k}, \la \rangle}.
$$
Using this notation, each of the aforementioned parts are
\begin{align*}
(\mbox{Structured part})&=w(\la)\frac{1}{m_d(B_R)}\sum_{p\in X_R}e^{-2\pi i\langle p,\lambda\rangle}\mu(p).\\
(\mbox{Almost Periodic part})&=\frac{1}{m_d(B_R)}\sum_{p\in X_R}e^{-2\pi i\langle p,\lambda\rangle}\mu(p) \, \Big(\E(\Xi_p(\la))-w(\la) \Big).\\
 (\mbox{Random part})  &=\frac{1}{m_d(B_R)}\sum_{p\in X_R}e^{-2\pi i\langle p,\lambda\rangle}\mu(p)\Big(\Xi_p(\la)-\E(\Xi_p(\la))\Big) \,.
\end{align*}

As we have already seen, by Proposition \ref{notPoisson}
\begin{equation*}
    \lim_{R\to\infty} \frac{w(\la)}{m_d(B_R)} \sum_{p\in X_R} e^{-2\pi i \langle p,\lambda \rangle}\mu(p) = \left\{ \begin{array}{lcc}
             w(\la)A(\lambda) &   if  & \lambda\in X^* \\
             \\ 0 &  if  & \lambda\notin X^*
             \end{array}
   \right.,
\end{equation*}
and by hypothesis \textbf{(\ref{eq:H5})} of asymptotic uniformity and the fact that $\mu$ is translation-bounded, we also have that
$$
   \lim_{R\to\infty}\frac{1}{m_d(B_R)}\sum_{p\in X_R}e^{-2\pi i\langle p,\lambda\rangle}\mu(p) \, \Big(\E(\Xi_p(\la))-w(\la) \Big)=0.
$$
So, it only remains to see that the random part also converges to zero as $R\to\infty$. The idea is to use Theorem \ref{chuecageneralizada}. So, let us check that the hypotheses are satisfied if we take $Y_p(\la)=\Xi_p(\la)-\e(\Xi_p(\la))$. Defined in this way, it is straightforward to verify that $\{Y_p(\lambda)\}$ satisfies the first two hypotheses of Theorem~\ref{chuecageneralizada}. Let us check only the last two of them, regarding the derivatives. Note that 
\begin{equation*}
    \frac{\partial Y_p(\la)}{\partial \la_j} = \e\left(2\pi i\sum_{k=1}^{n(p)} \alpha_{p,k} \,e^{-2\pi i \langle \xi_{p,k}, \la \rangle} (\xi_{p,k})_j\right) - 2\pi i\sum_{k=1}^{n(p)} \alpha_{p,k} \,e^{-2\pi i \langle \xi_{p,k}, \la \rangle} (\xi_{p,k})_j \,.
\end{equation*}
By Lemma \ref{aux1}, $\e\left(e^{-2\pi i \langle \xi_{p,k},\la \rangle} \xi_j\right)$ is continuous in $\la$. Moreover,
\begin{align*}
    \lVert \nabla Y_p(\la)\lVert_1 &\lesssim \e\left(\sum_{j=1}^d \sum_{k=1}^{n(p)} |\alpha_{p,k}||(\xi_{p,k})_j|  \right) + \sum_{j=1}^d \sum_{k=1}^{n(p)} |\alpha_{p,k}||(\xi_{p,k})_j| \\[2mm]
    &\lesssim\e\left(\sum_{k=1}^{n(p)} |\alpha_{p,k}|\| \xi_{p,k}\|\right) +  \sum_{k=1}^{n(p)} |\alpha_{p,k}|\| \xi_{p,k}\| \\
    &\lesssim \e(\| \xi\|) +  \sum_{k=1}^{n(p)} |\alpha_{p,k}|\| \xi_{p,k}\| \,.
\end{align*}
Since
\begin{equation*}
   \e\left(\left(\sum_{k=1}^{n(p)} |\alpha_{p,k}|\| \xi_{p,k}\|\right)^2\right) \,
\end{equation*}
is uniformly bounded (by \textbf{(\ref{eq:H6})}), the last condition on $\{Y_p\}$ is satisfied. Therefore, by Theorem \ref{chuecageneralizada}, the random part also converges to zero as $R\to \infty$.

\section*{Auxiliary Results}

\begin{appendix}

\section{An Ad-Hoc Law of Large Numbers}

Here we present ideas inspired by Theorem 6 in \cite{MR978305} tailored to our setting.

\begin{teor}\label{teoadoc}
Let \(\{X_n\}\) be a sequence of complex random variables such that $\e(|X_n|^2)$ is uniformly bounded over $n$, and let \(\{\alpha_n\}\) be a sequence of nonzero complex numbers. Suppose there exists an increasing sequence \(\{N_k\}\) such that
\[
\sum_{k\ge1}\mathbb{E}\left(\left|\frac{\sum_{n=1}^{N_k}\alpha_nX_n}{\sum_{n=1}^{N_k}|\alpha_n|}\right|^2\right)<\infty \,.
\]
If
\begin{equation*}
       \sum_{k\geq 1} \left(\frac{\sum_{n=1}^{N_{k+1}} |\a_n|}{\sum_{n=1}^{N_{k}} |\a_n|}-1\right)^2 < \infty \,,
    \end{equation*}
then almost surely
\[
\lim_{N\to\infty}\frac{\sum_{n=1}^{N}\alpha_nX_n}{\sum_{n=1}^{N}|\alpha_n|}=0 \,.
\]
\end{teor}

\begin{proof}
Since 
\[
\e\left(\sum_{k\ge1}\left|\frac{\sum_{n=1}^{N_k}\alpha_nX_n}{\sum_{n=1}^{N_k}|\alpha_n|}\right|^2\right) = \sum_{k\ge1}\mathbb{E}\left(\left|\frac{\sum_{n=1}^{N_k}\alpha_nX_n}{\sum_{n=1}^{N_k}|\alpha_n|}\right|^2\right)<\infty \,,
\]
then 
\[
\sum_{k\ge1}\left|\frac{\sum_{n=1}^{N_k}\alpha_nX_n}{\sum_{n=1}^{N_k}|\alpha_n|}\right|^2
\]
is almost surely finite, and so 
\[
\lim_{k\to\infty}\frac{\sum_{n=1}^{N_k}\alpha_nX_n}{\sum_{n=1}^{N_k}|\alpha_n|}=0 
\]
almost surely. Let
    \begin{equation*}
        \lVert X\rVert_2:=\sqrt{\e(|X|^2)} \,.
    \end{equation*}
We have
    \begin{align*}
        \left\lVert\max_{1\leq s< N_{k+1}-N_k}\left|\frac{\sum_{n=N_k+1}^{N_k+s} \a_n X_n}{\sum_{n=1}^{N_k} |\a_n|}\right|\right\rVert_2 &\leq \left\lVert\max_{1\leq s< N_{k+1}-N_k}\frac{\sum_{n=N_k+1}^{N_k+s} |\a_n| |X_n|}{\sum_{n=1}^{N_k} |\a_n|}\right\rVert_2  \\[2mm]
        &\leq \left\lVert\frac{\sum_{n=N_k+1}^{N_{k+1}} |\a_n| |X_n|}{\sum_{n=1}^{N_k} |\a_n|}\right\rVert_2 \\[2mm]
        &\lesssim\frac{\sum_{n=N_k+1}^{N_{k+1}} \left|\a_n\right| }{\sum_{n=1}^{N_k} |\a_n|} \\[2mm]
        &= \left(\frac{\sum_{n=1}^{N_{k+1}} |\a_n|}{\sum_{n=1}^{N_{k}} |\a_n|}-1\right)  \,.
    \end{align*}
Hence
    \begin{equation*}
        \sum_{k\geq 1} \e\left(\left(\max_{1\leq s< N_{k+1}-N_k}\left|\frac{\sum_{n=N_k+1}^{N_k+s} \a_n X_n}{\sum_{n=1}^{N_k} |\a_n|}\right|\right)^2\right) \lesssim \sum_{k\geq 1} \left(\frac{\sum_{n=1}^{N_{k+1}} |\a_n|}{\sum_{n=1}^{N_{k}} |\a_n|}-1\right)^2 < \infty \,.
    \end{equation*}
As before, 
\[
\lim_{k\to\infty}\max_{1\leq s< N_{k+1}-N_k}\left|\frac{\sum_{n=N_k+1}^{N_k+s} \a_n X_n}{\sum_{n=1}^{N_k} |\a_n|}\right|=0 
\]
almost surely. Finally, if $N_k\leq N < N_{k+1}$
    \begin{equation*}
        \left|\frac{\sum_{n=1}^{N} \a_n X_n}{\sum_{n=1}^{N} |\a_n|}\right| \leq \left|\frac{\sum_{n=1}^{N} \a_n X_n}{\sum_{n=1}^{N_k} |\a_n|}\right| \leq \left|\frac{\sum_{n=1}^{N_k} \a_n X_n}{\sum_{n=1}^{N_k} |\a_n|}\right| + \max_{1\leq s< N_{k+1}-N_k}\left|\frac{\sum_{n=N_k+1}^{N_k+s} \a_n X_n}{\sum_{n=1}^{N_k} |\a_n|}\right| \to 0 \,.
    \end{equation*}

\end{proof}

\begin{teor}\label{leyadoc}
Let \(\mu\) be a translation-bounded purely atomic measure with support \(X\subseteq\mathbb{R}^d\) locally finite such that
\[
|\mu|(X)=\infty \,.
\]
Let \(\{Y_p\}_{p\in X}\) be complex random variables with zero mean satisfying \(\mathbb{E}(|Y_p|^2)\) is uniformly bounded over $p\in X$. Suppose further that there exists \(R_0\) such that whenever \(\|p-p'\|_2>R_0\), the variables \(Y_p\) and \(Y_{p'}\) are independent. Then almost surely
\[
\lim_{R\to\infty}\frac{1}{|\mu|(B_R)}\sum_{p\in X_R}\mu(p)Y_p=0 \,
\]
and 
\[
\lim_{R\to\infty}\frac{1}{|\mu|(B_R)}\sum_{p\in X_R}|\mu(p)|Y_p=0 \,.
\]
\end{teor}

\begin{proof}
Enumerate the points of \(X\) (that are infinite) as \(\{p_n\}_{n\in\mathbb{N}}\) in non-decreasing order of \(\|p_n\|\). As $X$ is locally finite, for any complex sequence $\{Z_p\}_{p\in X}$ 
\[
\lim_{R\to\infty}\frac{1}{|\mu|(B_R)}\sum_{p\in X_R}Z_p=0 \, \iff \lim_{N\to\infty}\frac{1}{\sum_{n=1}^N |\mu(p_n)|}\sum_{n=1}^N Z_{p_n}=0 \,.
\]
Define
\[
N_k:=\min\left\{N\in\mathbb{N}:\sum_{n=1}^N|\mu(p_n)|\ge k^2\right\}.
\]
Since the \(\mu(p)\) are bounded, for sufficiently large \(k_0\) we have \(N_k<N_{k+1}\). Then 
\[
\frac{\sum_{n=1}^{N_{k+1}}|\mu(p_n)|}{\sum_{n=1}^{N_k}|\mu(p_n)|}\le \frac{(k+2)^2}{k^2} \,, \quad \text{ for all } k\geq k_0 \,.
\]
We have
    \begin{equation*}
        \sum_{k\geq k_0}\left(\frac{\sum_{n=1}^{N_{k+1}} |\mu(p_n)|}{\sum_{n=1}^{N_{k}} |\mu(p_n)|} - 1\right)^2 \leq \sum_{k\geq k_0} \left( \frac{(k+2)^2}{k^2}-1\right)^2 = 16\sum_{k\geq k_0} \frac{(k+1)^2}{k^4} < \infty \,.
    \end{equation*}
Also
     \begin{align*}
         &\var\left(\frac{1}{\sum_{n=1}^{N_k} |\mu(p_n)|}\sum_{n=1}^{N_k} \mu(p_n)Y_{p_n}\right) = \\
         &\frac{1}{\left(\sum_{n=1}^{N_k} |\mu(p_n)|\right)^2}\left(\sum_{i=1}^{N_k} \sum_{j=1}^{N_k} \mu(p_i)\overline{
         \mu(p_j)}\cov(Y_{p_i},Y_{p_j})\right) \leq\\
         & \frac{1}{\left(\sum_{n=1}^{N_k} |\mu(p_n)|\right)^2} \left(\sum_{i=1}^ {N_k}|\mu(p_i)|\sum_{j: \lVert p_j - p_i \rVert \leq R_0}|\mu(p_j)|\left|\cov(Y_{p_i},Y_{p_j})\right| \right) \underset{(\text{Cauchy-Schwarz})}{\lesssim}\\
         & \frac{1}{\left(\sum_{n=1}^{N_k} |\mu(p_n)|\right)^2} \left(\sum_{i=1}^ {N_k}|\mu(p_i)|\sum_{j: \lVert p_j - p_i \rVert \leq R_0}|\mu(p_j)|\right) \leq\\
         &\frac{1}{\left(\sum_{n=1}^{N_k} |\mu(p_n)|\right)^2} \left(\sum_{i=1}^{N_k}|\mu(p_i)||\mu|(B_{R_0}(p_i)) \right) \lesssim\\
         &\frac{1}{\sum_{n=1}^{N_k} |\mu(p_n)|} \,.
     \end{align*}
Apply Theorem \ref{teoadoc} and the first desired limit is obtained. The second limit follows the same proof.
\end{proof}

\section*{Acknowledgments}

The authors thank the two anonymous reviewers for their careful
reading of the manuscript and their valuable comments,
which helped improve the presentation of this work.

\end{appendix}

\printbibliography

@article {MR3939563,
    AUTHOR = {Lenz, Daniel and Strungaru, Nicolae},
     TITLE = {On weakly almost periodic measures},
   JOURNAL = {Trans. Amer. Math. Soc.},
  FJOURNAL = {Transactions of the American Mathematical Society},
    VOLUME = {371},
      YEAR = {2019},
    NUMBER = {10},
     PAGES = {6843--6881},
      ISSN = {0002-9947,1088-6850},
   MRCLASS = {43A05 (11K70 37B05 37N20 43A60)},
  MRNUMBER = {3939563},
MRREVIEWER = {Anatoly\ N.\ Kochubei},
       DOI = {10.1090/tran/7422},
       URL = {https://doi.org/10.1090/tran/7422},
}

@article{MR2876415,
	author = {Meyer, Yves},
	fjournal = {African Diaspora Journal of Mathematics},
	issn = {1539-854X},
	journal = {Afr. Diaspora J. Math.},
	mrclass = {52C23 (11K16 42A75)},
	mrnumber = {2876415},
	mrreviewer = {Anton\ Vladimirovich\ Shutov},
	number = {1},
	pages = {1--45},
	title = {Quasicrystals, almost periodic patterns, mean-periodic functions and irregular sampling},
	volume = {13},
	year = {2012}}

@article{MR4866154,
	author = {Kulikov, Aleksei and Nazarov, Fedor and Sodin, Mikhail},
	fjournal = {Journal of Mathematical Physics, Analysis, Geometry},
	issn = {1812-9471,1817-5805},
	journal = {J. Math. Phys. Anal. Geom.},
	mrclass = {42A38 (42B10)},
	mrnumber = {4866154},
	mrreviewer = {Debkumar\ Giri},
	number = {1},
	pages = {84--130},
	title = {Fourier uniqueness and non-uniqueness pairs},
	volume = {21},
	year = {2025}}

@article{Fa24,
	author = {Favorov, S. Yu.},
	doi = {10.1007/s10476-024-00031-y},
	issn = {1588-273X},
	journal = {Analysis Mathematica},
	month = jun,
	number = {2},
	pages = {455--462},
	publisher = {Springer Science and Business Media LLC},
	title = {A crystalline measure that is not a Fourier quasicrystal},
	url = {http://dx.doi.org/10.1007/s10476-024-00031-y},
	volume = {50},
	year = {2024}}

@incollection{MR4939386,
	author = {Meyer, Yves},
	booktitle = {{\it {S}electa}. {Y}ves {M}eyer. {V}ol. 1},
	isbn = {978-2-85629-992-0},
	mrclass = {42A75 (11M06 30B50)},
	mrnumber = {4939386},
	pages = {131--140},
	publisher = {Soc. Math. France, Paris},
	series = {Doc. Math. (Paris)},
	title = {Zeta functions and crystalline measures},
	volume = {22},
	year = {[2025] \copyright 2025}}

@misc{boyvalenkov2025growthmassescrystallinemeasures,
	archiveprefix = {arXiv},
	author = {Peter Boyvalenkov and Sergii Yu. Favorov},
	eprint = {2503.19567},
	primaryclass = {math.FA},
	title = {Growth of masses of crystalline measures},
	url = {https://arxiv.org/abs/2503.19567},
	year = {2025}}

@misc{ArXiv-2312.11185,
	archiveprefix = {arXiv},
	author = {Felipe Gon{\c c}alves},
	eprint = {2312.11185},
	primaryclass = {math.CA},
	title = {A Classification of Fourier Summation Formulas and Crystalline Measures},
	url = {https://arxiv.org/abs/2312.11185},
	year = {2023}}

@book{MR1995517,
	author = {Assani, Idris},
	doi = {10.1142/4538},
	isbn = {981-02-4439-8},
	mrclass = {37A30 (28D05 47A35 82B03)},
	mrnumber = {1995517},
	mrreviewer = {U. Krengel},
	pages = {xii+216},
	publisher = {World Scientific Publishing Co. Inc. River Edge, NJ},
	title = {Wiener Wintner ergodic theorems},
	url = {https://doi.org/10.1142/4538},
	year = {2003}}

@article{MR472749,
	author = {Coquet, Jean and Kamae, Teturo and Mend{\`e}s France, Michel},
	fjournal = {Bulletin de la Soci\'{e}t\'{e} Math\'{e}matique de France},
	issn = {0037-9484},
	journal = {Bull. Soc. Math. France},
	mrclass = {10K35 (42A84)},
	mrnumber = {472749},
	mrreviewer = {J.-P. Kahane},
	number = {4},
	pages = {369--384},
	title = {Sur la mesure spectrale de certaines suites arithm\'{e}tiques},
	url = {http://www.numdam.org/item?id=BSMF_1977__105__369_0},
	volume = {105},
	year = {1977}}

@book{MR4697523,
	author = {Kurasov, Pavel},
	doi = {10.1007/978-3-662-67872-5},
	isbn = {978-3-662-67870-1; 978-3-662-67872-5},
	mrclass = {81-02 (05-02 34B45 47A10 81Q35 81U40)},
	mrnumber = {4697523},
	mrreviewer = {Mohammed El A\"{\i}di, Universidad Nacional de Colombia},
	pages = {xvi+639},
	publisher = {Birkh\"{a}user/Springer, Berlin},
	series = {Operator Theory: Advances and Applications},
	title = {Spectral geometry of graphs},
	url = {https://doi.org/10.1007/978-3-662-67872-5},
	volume = {293},
	year = {2024}}

@article{MR4576026,
	author = {Meyer, Yves},
	fjournal = {Ast\'{e}risque},
	isbn = {978-2-85629-968-5},
	issn = {0303-1179},
	journal = {Ast\'{e}risque},
	mrclass = {42-02 (11H31 42A38 46F10 52C23)},
	mrnumber = {4576026},
	mrreviewer = {Joseph D. Lakey},
	number = {438, S\'{e}minaire Bourbaki. Vol. 2021/2022. Expos\'{e}s 1181--1196},
	pages = {Exp. No. 1194, 479--494},
	title = {Mesures cristallines et applications [{\it d'apr{\`e}s} {P}. {K}urasov, {N}. {L}ev, {A}. {O}levskii, {P}. {S}arnak, et {M}. {V}iazovska]},
	url = {https://mathscinet.ams.org/mathscinet-getitem?mr=4576026},
	year = {2022}}

@article{MR4841781,
	author = {Alon, Lior and Kummer, Mario and Kurasov, Pavel and Vinzant, Cynthia},
	doi = {10.1007/s00222-024-01307-8},
	fjournal = {Inventiones Mathematicae},
	issn = {0020-9910},
	journal = {Invent. Math.},
	mrclass = {52C23 (14P05 42A38 42A75 42B10)},
	mrnumber = {4841781},
	number = {1},
	pages = {321--376},
	title = {Higher dimensional {F}ourier quasicrystals from {L}ee-{Y}ang varieties},
	url = {https://doi.org/10.1007/s00222-024-01307-8},
	volume = {239},
	year = {2025}}

@article{MR4678143,
	author = {Ramos, Jo\~{a}o P. G. and Sousa, Mateus},
	doi = {10.2140/apde.2023.16.2327},
	fjournal = {Analysis \& PDE},
	issn = {2157-5045},
	journal = {Anal. PDE},
	mrclass = {42A38 (11F03 41A05 46E39)},
	mrnumber = {4678143},
	mrreviewer = {Stefan Steinerberger},
	number = {10},
	pages = {2327--2384},
	title = {Perturbed interpolation formulae and applications},
	url = {https://doi.org/10.2140/apde.2023.16.2327},
	volume = {16},
	year = {2023}}

@book{MR1157815,
	author = {Rudin, Walter},
	edition = {Second},
	isbn = {0-07-054236-8},
	mrclass = {46-01 (47-01)},
	mrnumber = {1157815},
	pages = {xviii+424},
	publisher = {McGraw-Hill, Inc., New York},
	series = {International Series in Pure and Applied Mathematics},
	title = {Functional analysis},
	url = {https://mathscinet.ams.org/mathscinet-getitem?mr=1157815},
	year = {1991}}

@article{MR4406130,
	author = {Yakir, Oren},
	doi = {10.1093/imrn/rnaa316},
	fjournal = {International Mathematics Research Notices. IMRN},
	issn = {1073-7928},
	journal = {Int. Math. Res. Not. IMRN},
	mrclass = {60D05 (11H06 11P21 52C07)},
	mrnumber = {4406130},
	number = {8},
	pages = {6243--6261},
	title = {Recovering the lattice from its random perturbations},
	url = {https://doi.org/10.1093/imrn/rnaa316},
	year = {2022}}

@incollection{MR1798991,
	author = {Schlottmann, Martin},
	booktitle = {Directions in mathematical quasicrystals},
	mrclass = {52C23 (11K70 37A05 43A25 82D25)},
	mrnumber = {1798991},
	mrreviewer = {Richard Kenyon},
	pages = {143--159},
	publisher = {Amer. Math. Soc., Providence, RI},
	series = {CRM Monogr. Ser.},
	title = {Generalized model sets and dynamical systems},
	url = {https://mathscinet.ams.org/mathscinet-getitem?mr=1798991},
	volume = {13},
	year = {2000}}

@article{MR4526342,
	author = {Petrache, Mircea and Viera, Rodolfo},
	doi = {10.1007/s10955-022-03059-2},
	fjournal = {Journal of Statistical Physics},
	issn = {0022-4715},
	journal = {J. Stat. Phys.},
	mrclass = {60D05 (11H06 52C23)},
	mrnumber = {4526342},
	number = {2},
	pages = {Paper No. 39, 26},
	title = {Almost sure recovery in quasi-periodic structures},
	url = {https://doi.org/10.1007/s10955-022-03059-2},
	volume = {190},
	year = {2023}}

@incollection{MR3791850,
	author = {Moody, Robert V. and Strungaru, Nicolae},
	booktitle = {Aperiodic order. {V}ol. 2},
	mrclass = {42A75 (28A12 28C15)},
	mrnumber = {3791850},
	pages = {173--270},
	publisher = {Cambridge Univ. Press, Cambridge},
	series = {Encyclopedia Math. Appl.},
	title = {Almost periodic measures and their {F}ourier transforms},
	url = {https://mathscinet.ams.org/mathscinet-getitem?mr=3791850},
	volume = {166},
	year = {2017}}

@incollection{MR4633114,
	author = {Meyer, Yves},
	booktitle = {Theoretical physics, wavelets, analysis, genomics---an indisciplinary tribute to {A}lex {G}rossmann},
	doi = {10.1007/978-3-030-45847-8\_17},
	mrclass = {42A32 (43A46)},
	mrnumber = {4633114},
	pages = {389--407},
	publisher = {Birkh\"{a}user/Springer, Cham},
	series = {Appl. Numer. Harmon. Anal.},
	title = {Curved model sets and crystalline measures},
	url = {https://doi.org/10.1007/978-3-030-45847-8_17},
	year = {[2023] \copyright 2023}}

@article{MR978305,
	author = {Lyons, Russell},
	doi = {10.1307/mmj/1029003816},
	fjournal = {Michigan Mathematical Journal},
	issn = {0026-2285},
	journal = {Michigan Math. J.},
	mrclass = {60F15},
	mrnumber = {978305},
	mrreviewer = {Andrew D. Barbour},
	number = {3},
	pages = {353--359},
	title = {Strong laws of large numbers for weakly correlated random variables},
	url = {https://doi.org/10.1307/mmj/1029003816},
	volume = {35},
	year = {1988}}

@article{MR106208,
	author = {Kahane, J.-P. and Mandelbrojt, S.},
	fjournal = {Annales Scientifiques de l'\'{E}cole Normale Sup\'{e}rieure. Troisi{\`e}me S\'{e}rie},
	issn = {0012-9593},
	journal = {Ann. Sci. \'{E}cole Norm. Sup. (3)},
	mrclass = {10.00},
	mrnumber = {106208},
	mrreviewer = {K. Chandrasekharan},
	pages = {57--80},
	title = {Sur l'\'{e}quation fonctionnelle de {R}iemann et la formule sommatoire de {P}oisson},
	url = {http://www.numdam.org/item?id=ASENS_1958_3_75_1_57_0},
	volume = {75},
	year = {1958}}

@article{MR2483565,
	author = {Dyson, Freeman},
	fjournal = {Notices of the American Mathematical Society},
	issn = {0002-9920},
	journal = {Notices Amer. Math. Soc.},
	mrclass = {01A70},
	mrnumber = {2483565},
	number = {2},
	pages = {212--223},
	title = {Birds and frogs},
	url = {https://mathscinet.ams.org/mathscinet-getitem?mr=2483565},
	volume = {56},
	year = {2009}}

@article{MR3338010,
	author = {Lev, Nir and Olevskii, Alexander},
	doi = {10.1007/s00222-014-0542-z},
	fjournal = {Inventiones Mathematicae},
	issn = {0020-9910},
	journal = {Invent. Math.},
	mrclass = {28A12 (28A05)},
	mrnumber = {3338010},
	mrreviewer = {Omer Edhan},
	number = {2},
	pages = {585--606},
	title = {Quasicrystals and {P}oisson's summation formula},
	url = {https://doi.org/10.1007/s00222-014-0542-z},
	volume = {200},
	year = {2015}}

@article{MR3448912,
	author = {Kolountzakis, Mihail N.},
	doi = {10.1007/s00041-015-9416-z},
	fjournal = {The Journal of Fourier Analysis and Applications},
	issn = {1069-5869},
	journal = {J. Fourier Anal. Appl.},
	mrclass = {42A38},
	mrnumber = {3448912},
	mrreviewer = {Tian-Xiao He},
	number = {1},
	pages = {1--5},
	title = {Fourier pairs of discrete support with little structure},
	url = {https://doi.org/10.1007/s00041-015-9416-z},
	volume = {22},
	year = {2016}}

@article{MR1328260,
	author = {Hof, A.},
	fjournal = {Communications in Mathematical Physics},
	issn = {0010-3616},
	journal = {Comm. Math. Phys.},
	mrclass = {82D25 (52C22 78A45)},
	mrnumber = {1328260},
	mrreviewer = {Tatiana Allen},
	number = {1},
	pages = {25--43},
	title = {On diffraction by aperiodic structures},
	url = {http://projecteuclid.org/euclid.cmp/1104272610},
	volume = {169},
	year = {1995}}

@article{MR3815253,
	author = {Torquato, Salvatore},
	doi = {10.1016/j.physrep.2018.03.001},
	fjournal = {Physics Reports. A Review Section of Physics Letters},
	issn = {0370-1573},
	journal = {Phys. Rep.},
	mrclass = {82B44},
	mrnumber = {3815253},
	pages = {1--95},
	title = {Hyperuniform states of matter},
	url = {https://doi.org/10.1016/j.physrep.2018.03.001},
	volume = {745},
	year = {2018}}

@article{MR4735950,
	author = {Bj\"{o}rklund, Michael and Hartnick, Tobias},
	doi = {10.1007/s00208-023-02647-1},
	fjournal = {Mathematische Annalen},
	issn = {0025-5831},
	journal = {Math. Ann.},
	mrclass = {82C44 (11K06 43A25 60G55)},
	mrnumber = {4735950},
	mrreviewer = {C\'{e}sar R. de Oliveira},
	number = {1},
	pages = {365--426},
	title = {Hyperuniformity and non-hyperuniformity of quasicrystals},
	url = {https://doi.org/10.1007/s00208-023-02647-1},
	volume = {389},
	year = {2024}}

@article{MR3949027,
	author = {Radchenko, Danylo and Viazovska, Maryna},
	doi = {10.1007/s10240-018-0101-z},
	fjournal = {Publications Math\'{e}matiques. Institut de Hautes \'{E}tudes Scientifiques},
	issn = {0073-8301},
	journal = {Publ. Math. Inst. Hautes \'{E}tudes Sci.},
	mrclass = {11F37 (11F03 11F11 42A10)},
	mrnumber = {3949027},
	mrreviewer = {Jeremy Lovejoy},
	pages = {51--81},
	title = {Fourier interpolation on the real line},
	url = {https://doi.org/10.1007/s10240-018-0101-z},
	volume = {129},
	year = {2019}}

@article{MR4294062,
	author = {Bakan, Andrew and Hedenmalm, Haakan and Montes-Rodr\'{\i}guez, Alfonso and Radchenko, Danylo and Viazovska, Maryna},
	doi = {10.1073/pnas.2023227118},
	fjournal = {Proceedings of the National Academy of Sciences of the United States of America},
	issn = {0027-8424},
	journal = {Proc. Natl. Acad. Sci. USA},
	mrclass = {42B37},
	mrnumber = {4294062},
	number = {15},
	pages = {Paper No. 2023227118, 4},
	title = {Fourier uniqueness in even dimensions},
	url = {https://doi.org/10.1073/pnas.2023227118},
	volume = {118},
	year = {2021}}

@book{MR3136260,
	author = {Baake, Michael and Grimm, Uwe},
	doi = {10.1017/CBO9781139025256},
	isbn = {978-0-521-86991-1},
	mrclass = {52C23 (11H06 20C35 20H15 82D25)},
	mrnumber = {3136260},
	mrreviewer = {Jean-Pierre Gazeau},
	note = {A mathematical invitation, With a foreword by Roger Penrose},
	pages = {xvi+531},
	publisher = {Cambridge University Press, Cambridge},
	series = {Encyclopedia of Mathematics and its Applications},
	title = {Aperiodic order. {V}ol. 1},
	url = {https://doi.org/10.1017/CBO9781139025256},
	volume = {149},
	year = {2013}}

@article{MR3482845,
	author = {Meyer, Yves},
	doi = {10.1073/pnas.1600685113},
	fjournal = {Proceedings of the National Academy of Sciences of the United States of America},
	issn = {0027-8424},
	journal = {Proc. Natl. Acad. Sci. USA},
	mrclass = {28A12},
	mrnumber = {3482845},
	mrreviewer = {Peter R. Massopust},
	number = {12},
	pages = {3152--3158},
	title = {Measures with locally finite support and spectrum},
	url = {https://doi.org/10.1073/pnas.1600685113},
	volume = {113},
	year = {2016}}

@article{MR4129870,
	author = {Kurasov, P. and Sarnak, P.},
	doi = {10.1063/5.0012286},
	fjournal = {Journal of Mathematical Physics},
	issn = {0022-2488},
	journal = {J. Math. Phys.},
	mrclass = {46F12 (42C05 46F05 46N50 52C23)},
	mrnumber = {4129870},
	mrreviewer = {S. K. Upadhyay},
	number = {8},
	pages = {083501, 13},
	title = {Stable polynomials and crystalline measures},
	url = {https://doi.org/10.1063/5.0012286},
	volume = {61},
	year = {2020}}

@article{MR4206541,
	author = {Olevskii, Alexander and Ulanovskii, Alexander},
	doi = {10.5802/crmath.142},
	fjournal = {Comptes Rendus Math\'{e}matique. Acad\'{e}mie des Sciences. Paris},
	issn = {1631-073X},
	journal = {C. R. Math. Acad. Sci. Paris},
	mrclass = {46F12 (42B10 52C23)},
	mrnumber = {4206541},
	number = {11-12},
	pages = {1207--1211},
	title = {Fourier quasicrystals with unit masses},
	url = {https://doi.org/10.5802/crmath.142},
	volume = {358},
	year = {2020}}

@article{MR4664429,
	author = {Alon, Lior and Cohen, Alex and Vinzant, Cynthia},
	doi = {10.1016/j.jfa.2023.110226},
	fjournal = {Journal of Functional Analysis},
	issn = {0022-1236},
	journal = {J. Funct. Anal.},
	mrclass = {32A08 (14T20 32A60 42C05)},
	mrnumber = {4664429},
	mrreviewer = {Sergiy U. Favorov},
	number = {2},
	pages = {Paper No. 110226, 10},
	title = {Every real-rooted exponential polynomial is the restriction of a {L}ee-{Y}ang polynomial},
	url = {https://doi.org/10.1016/j.jfa.2023.110226},
	volume = {286},
	year = {2024}}

\end{document}